\documentclass[leqno,11pt]{amsart}
\usepackage{amsmath,amscd}
\usepackage{epsfig,graphics}
\usepackage{amssymb,latexsym}
\usepackage{mathrsfs}
\usepackage{eucal}

\newtheorem{thm}{\bf Theorem}[section]
\newtheorem{df}[thm]{\bf Definition}
\newtheorem{prop}[thm]{\bf Proposition}
\newtheorem{cor}[thm]{\bf Corollary}

\newtheorem{rem}[thm]{\bf Remark}
\newtheorem{ex}[thm]{\bf Example}

\newcommand{\bs}{\boldsymbol}
\newcommand{\A}{\mathcal{A}}
\newcommand{\B}{\mathbf{B}}

\newcommand{\W}{\mathcal{W}}
\newcommand{\cP}{\mathscr{P}}

\newcommand{\pf}{\noindent{\bfseries Proof. }}

\newcommand{\w}{{\bf w}}
\newcommand{\bi}{\bs{\rm i}}
\newcommand{\bj}{\bs{\rm j}}

\newcommand{\bk}{\bs{\rm k}}
\newcommand{\bl}{\bs{\rm l}}

\newcommand{\M}{{\mathcal{M}}}
\newcommand{\F}{\mathcal{F}}

\newcommand{\N}{\mathcal{T}}
\newcommand{\gl}{\mathfrak{gl}}
\newcommand{\Z}{\mathbb{Z}}
\newcommand{\C}{\mathbb{C}}
\newcommand{\h}{\mathfrak{h}}
\newcommand{\te}{\widetilde{e}}
\newcommand{\tf}{\widetilde{f}}

\newcommand{\td}{\widetilde}

\numberwithin{equation}{section}

\begin{document}
\title[Crystal graph of a Generalized Verma module]
{Demazure crystals of generalized Verma modules and a flagged RSK
correspondence}
\author{JAE-HOON KWON}
\address{Department of Mathematics \\ University of Seoul   \\ Seoul 130-743, Korea }
\email{jhkwon@uos.ac.kr }

\thanks{This work was supported by KRF Grant 2008-314-C00004.
}

\begin{abstract}
We prove that the Robinson-Schensted-Knuth correspondence is a $\gl_{\infty}$-crystal
isomorphism between two realizations of the crystal graph of a
generalized Verma module with respect to a maximal
parabolic subalgebra of $\gl_{\infty}$.
A flagged version of the RSK correspondence
is derived in a natural way by computing a Demazure
crystal graph of a generalized Verma module. As an application, we discuss a relation between
a Demazure crystal  and plane partitions with a bounded condition.
\end{abstract} \maketitle

\section{Introduction}

The Robinson-Schensted-Knuth (simply RSK) algorithm \cite{Kn}  has
been playing fundamental roles in combinatorics and representation
theory with generalizations in various directions. It gives a
bijection between  the set $\M$ of $\mathbb{N}\times\mathbb{N}$
matrices with non-negative integers of finite support and  the set
$\N$ of pairs of semistandard tableaux of the same shape, and
explains in a bijective way the expansion of the Cauchy kernel into
Schur functions, called the Cauchy identity;
\begin{equation*}
\frac{1}{\prod_{i,j\geq
1}(1-x_iy_j)}=\sum_{\lambda}s_{\lambda}(X)s_{\lambda}(Y),
\end{equation*}
where the sum is over all partitions $\lambda$ and $s_{\lambda}(X)$
(or $s_{\lambda}(Y)$) is the Schur function in
$X=\{\,x_1,x_2,\ldots\,\}$ (or $Y=\{\,y_1,y_2,\ldots\,\}$). A
representation theoretic interpretation of the Cauchy identity can
be given by a general principle called Howe duality \cite{H}, that
is, $S(\mathbb{C}^{>0}\otimes \mathbb{C}^{>0})$, the symmetric
algebra generated by $\mathbb{C}^{>0}\otimes \mathbb{C}^{>0}$ has a
multiplicity-free decomposition into irreducible
$(\frak{gl}_{>0},\frak{gl}_{>0})$-bimodules parameterized by
partitions, where $\mathbb{C}^{>0}$ is the complex  vector space
with a basis $\{\,v_i\,|\,i\in\mathbb{N}\,\}$ and
$\frak{gl}_{>0}=\frak{gl}(\mathbb{C}^{>0})$ is the corresponding
general linear Lie algebra.

We have a more direct interpretation of the RSK map with the help of
the Kashiwara's crystal base theory of the quantum group
$U_q(\frak{gl}_{>0})$ \cite{Kas90,Kas94,KN}. That is,
both $\M$ and $\N$ have two $\frak{gl}_{>0}$-crystal structures
commuting with each other, which are called
$(\frak{gl}_{>0},\frak{gl}_{>0})$-bicrystals or double crystals
\cite{La}, and the RSK map is an isomorphism of bicrystals. The
decomposition as a $(\frak{gl}_{>0},\frak{gl}_{>0})$-bimodule
follows immediately by considering highest weight crystal elements
in $\M$. We refer the readers to \cite{DK,K07,La,Lc,Ln,vL} for more
results on bicrystal, its application and generalization to other
types of Lie algebras.

The main purpose of this paper is to give a new representation theoretic
interpretation of the
RSK correspondence and its applications.
Let $\frak{gl}_{\infty}$ be the general linear
Lie algebra, which is spanned by the elementary matrices $E_{ij}$
($i,j\in\mathbb{Z}\setminus\{0\}$). Let $\frak{l}=\frak{gl}_{<0}\oplus
\frak{gl}_{>0}$ be a Levi subalgebra, where $\frak{gl}_{\gtrless 0}$
is a subalgebra spanned by $E_{ij}$ ($i,j\gtrless0$), respectively,
and let $\frak{u}_\pm$ be the nilradical spanned by $E_{ij}$ for
$i>0, j<0$ (resp. $i<0, j>0$).
We may identify $S(\mathbb{C}^{>0}\otimes \mathbb{C}^{>0})$ with the
enveloping algebra $U(\frak{u}_-)$, which is a generalized Verma
module induced from the maximal parabolic subalgebra
$\frak{p}=\frak{l}\oplus \frak{u}_+$ (see \cite{Br} for a quantized
version of this fact and its relation with canonical basis).
Motivated by this observation, we introduce
$\frak{gl}_{\infty}$-crystal structures on $\M$ and $\N$ extending
the $(\frak{gl}_{<0}, \frak{gl}_{>0})$-bicrystal structures (note
that $\frak{gl}_{<0}\simeq \frak{gl}_{>0}$), and then show that the
RSK map is an isomorphism of $\frak{gl}_{\infty}$-crystals (Theorem
\ref{RSK iso}). Indeed, these are obtained by finding the missing
Kashiwara operators compatible with the RSK map, which correspond to
the simple root $\alpha_0$ connecting the Dynkin diagrams of
$\frak{gl}_{<0}$ and $\frak{gl}_{>0}$.

The RSK map also enables us to define a natural embedding of
$\B(n\Lambda_0)$ into $\M$, where $\Lambda_0$ is the $0$th
fundamental weight of $\gl_{\infty}$ and $\B(n\Lambda_0)$ is the
crystal graph of the irreducible $U_q(\frak{gl}_{\infty})$-module
with highest weight $n\Lambda_0$ (Proposition \ref{embedding1}).
Hence we may regard $\M$ as a crystal graph of the quantum group
$U_q(\frak{u}_-)$ since it can be realized as  a limit of
$\B(n\Lambda_0)$. In general, we give a combinatorial description of
a crystal graph of a generalized Verma module with arbitrary
$\frak{l}$-dominant highest weight.

Next, we define Demazure crystals $\M_w$ and $\N_w$ for $w\in W$
following \cite{Kas93}, where $W$ is the Weyl group for
$\frak{gl}_{\infty}$, and give explicit combinatorial descriptions
of them (Theorem \ref{Demazure-0} and \ref{Demazure-1}). As an
interesting corollary, the resulting flagged RSK correspondence
between $\M_w$ and $\N_w$ (Corollary \ref{main theorem}) explains a
nice relation between the support of a matrix in $\M$ and the flag
conditions of the corresponding tableaux in $\N$, which was observed
earlier in a purely combinatorial way (cf.\cite{St}). In terms of
characters, this can be summarized as follows; for each $w\in W$ we
have
$$\sum_{\substack{S\subset\mathbb{N}^2 \\
w(S)\leq w}}\prod_{(i,j)\in S}\frac{x_iy_j}{(1-x_iy_j)}=\sum_{\substack{\nu\in\cP \\
\ell(\nu)\leq
d(w)}}\widehat{s}_{\nu}(X_{\alpha(w)})\widehat{s}_{\nu}(Y_{\beta(w)}),$$
where the left hand side is the sum over supports in $\M$ dominated
by $w$ with respect to the Bruhat order, and the right-hand side is
the sum over products of flagged Schur functions with flag
conditions $\alpha(w),\beta(w)$ determined by $w$ (see Section 5 for
the precise definitions of these notations). We present variations
by considering symmetric matrices in $\M$ as crystal graphs for
affine Lie subalgebras $\frak{b}_{\infty}$ and $\frak{c}_{\infty}$
of $\frak{gl}_{\infty}$.

Finally, we discuss an application to plane partitions. We show that
a Demazure crystal associated with $w$ corresponds to a set of
(symmetric) plane partitions whose shapes are bounded by a partition
$\lambda$ corresponding to $w$, and obtain its trace generating
function as Demazure characters.

The paper is organized as follows. In Section 2, we recall necessary
background on crystal graphs. In Section 3, we define
$\frak{gl}_{\infty}$-crystal structures on $\M$ and $\N$. In Section
4, we give a combinatorial description of a crystal graph of a
generalized Verma module with arbitrary $\frak{l}$-dominant highest
weights including $\M$. In Section 5, we define and compute the
Demazure crystals $\M_w$ an $\N_w$ explicitly. In Section 6, we
discuss an application of Demazure crystals to plane
partitions.\vskip 3mm

{\bf Acknowledgement} The author would like to thank A. Lascoux for
helpful comments  and discussion on this work.

\section{Preliminaries}
\subsection{Lie algebra $\gl_{\infty}$}
Let $\Z^{\times}$ denote the set of non-zero integers. Let
$\gl_\infty$ denote the Lie algebra of $\Z^{\times}\times \Z^{\times}$ complex matrices
with finitely many non-zero entries.
Let $E_{ij}$ be the elementary matrix with $1$ at the $i$-th row and
the $j$-th column and zero elsewhere.

The Cartan subalgebra is given by $\h=\bigoplus_{i\in \Z^{\times}}\C
E_{ii}$. Let $\Pi^{\vee}=\{\, h_{-i}=E_{-i-1, -i-1}-E_{-i, -i},\
h_i=E_{ii}-E_{i+1,i+1}\, (i\in \Z_{>0}),\ h_0=E_{-1-1}-E_{11} \, \}$
be the set of simple coroots, $\Pi=\{\,
\alpha_{-i}=\epsilon_{-i-1}-\epsilon_{-i}\,\
\alpha_i=\epsilon_i-\epsilon_{i+1} \ (i\in \Z_{>0}),
\alpha_0=\epsilon_{-1}-\epsilon_1\, \}$ the set of simple roots, and
$\Delta^+ =\{\, \epsilon_i-\epsilon_j\,|\, i,j\in\Z^{\times}, i<j
\,\}$ the set of positive roots, where $\epsilon_i\in\h^*$ is
determined by $\langle \epsilon_i,E_{jj}\rangle=\delta_{ij}$ and
$\langle\cdot,\cdot\rangle$ is a natural pairing on
$\frak{h}^*\times\frak{h}$. The Dynkin diagram associated with the
Cartan matrix
$\left(\langle\alpha_j,h_i\rangle\right)_{i,j\in\mathbb{Z}}$ is
\begin{center}
\hskip 2cm\setlength{\unitlength}{0.45cm}
\begin{picture}(15,2)(0,0)
\put(-.85,.5){\line(1,0){1}}
\put(0,0){\makebox(1,1){$\bigcirc$}}\put(.85,.5){\line(1,0){1}}
\put(5,0){\makebox(1,1){$\bigcirc$}}\put(3.85,.5){\line(1,0){1.2}}\put(5.85,.5){\line(1,0){1.2}}
\put(7,0){\makebox(1,1){$\bigcirc$}}\put(7.85,.5){\line(1,0){1.2}}
\put(9,0){\makebox(1,1){$\bigcirc$}}\put(9.85,.5){\line(1,0){1.2}}
\put(14,0){\makebox(1,1){$\bigcirc$}}\put(12.85,.5){\line(1,0){1.2}}\put(14.85,.5){\line(1,0){1.2}}
\put(2.5,.2){$\cdots$}\put(11.5,.2){$\cdots$}
\put(-2,.2){$\cdots$}\put(16.5,.2){$\cdots$}

\put(0.5,-.5){\makebox(0,0)[c]{\tiny ${\alpha}_{-n}$}}
\put(5.5,-.5){\makebox(0,0)[c]{\tiny ${\alpha}_{-1}$}}
\put(7.5,-.5){\makebox(0,0)[c]{\tiny ${\alpha}_0$}}
\put(9.5,-.5){\makebox(0,0)[c]{\tiny ${\alpha}_1$}}
\put(14.5,-.5){\makebox(0,0)[c]{\tiny ${\alpha}_n$}}

\end{picture}\hskip 1.2cm .
\end{center}\vskip 2mm

Let $Q=\bigoplus_{i\in\mathbb{Z}}\mathbb{Z}\alpha_i$ be the
root lattice. Let
$P=\bigoplus_{i\in\Z^{\times}}\Z\epsilon_i\oplus\Z\Lambda_0$ be the
weight lattice of $\gl_{\infty}$, where $\Lambda_0$ is given by
$\langle\Lambda_0,E_{-i,-i}\rangle=-\langle\Lambda_0,E_{i,i}\rangle=\frac{1}{2}$
for $i\in\Z_{>0}$. There is a partial ordering $\geq$ on
$P$, where $\lambda\geq \mu $ if and only if $\lambda-\mu$ is a
non-negative integral linear combination of $\alpha_i$'s ($i\in\Z$).
Let $P^+=\{\,\Lambda\in P\,|\,\Lambda(h_i)\geq 0,\ i\in\Z\,\}$,
the set of dominant integral weights. For $i\in \Z^{\times}$, let
$\Lambda_i$ be
\begin{align*}
\Lambda_i=
\begin{cases}
\Lambda_0-\sum_{k=i}^{-1}\epsilon_k, & \text{if  $i< 0$}, \\
\Lambda_0+\sum_{k=1}^i\epsilon_k, & \text{if  $i> 0$},
\end{cases}
\end{align*}
We call $\Lambda_i\in P^+$ ($i\in\mathbb{Z}$) the $i$-th fundamental weight of
${\gl}_\infty$.

\subsection{Crystal graphs}
Let us recall the notion of crystal graphs for the Lie algebra
$\gl_{\infty}$ (cf. \cite{Kas90,Kas93,Kas94}).
\begin{df}\label{crystal graph}\mbox{}
{\rm  A  {\it $\gl_{\infty}$-crystal} is a set $B$ together with the
maps ${\rm wt}  : B \rightarrow P$, $\varepsilon_i, \varphi_i: B
\rightarrow \mathbb{Z}\cup \{-\infty\}$ and $\te_i, \tf_i: B
\rightarrow B\cup\{{\bf 0}\}$ for $i\in \Z$, satisfying the
following conditions;
\begin{itemize}
\item[(1)] for $b\in B$, we have
$$\varphi_i(b) =\langle {\rm wt}(b),h_i \rangle +
\varepsilon_i(b),$$

\item[(2)] if $\te_i b \in B$ for $b\in B$, then
$$\varepsilon_i(\te_i b) = \varepsilon_i(b) - 1,\ \varphi_i(\te_i b) =
\varphi_i(b) + 1,\ {\rm wt}(\te_ib)={\rm wt}(b)+\alpha_i,$$

\item[(3)] if $\tf_i b \in B$ for $b\in B$, then
$$\varepsilon_i(\tf_i b) = \varepsilon_i(b) + 1,\ \varphi_i(\tf_i b) =
\varphi_i(b) - 1,\ {\rm wt}({\tf_i}b)={\rm wt}(b)-\alpha_i,$$

\item[(4)] $\tf_i b = b'$ if and only if $b = \te_i
b'$ for all $i\in \Z$, $b, b' \in B$,

\item[(5)] If $\varphi_i(b)=-\infty$, then $\te_ib=\tf_ib={\bf
0}$,
\end{itemize}
where ${\bf 0}$ is a formal symbol and $-\infty$ is the smallest
element in $\mathbb{Z}\cup\{-\infty\}$ such that $-\infty+n=-\infty$
for all $n\in\mathbb{Z}$.}
\end{df}

A $\gl_{\infty}$-crystal $B$ becomes a colored oriented graph, where
$b\stackrel{i}{\rightarrow}b'$ if and only if $b'=\tf_{i}b$ $(i\in
\Z)$, and it is called a {\it crystal graph for $\gl_{\infty}$}. Let
$\mathbb{C}[P]$ be the group algebra of $P$ with basis
$\{\,e^{\lambda}\,|\,\lambda\in P\,\}$. We define the {\it character
of $B$} by ${\rm ch}B=\sum_{b\in B}e^{{\rm wt}(b)}$.

Let $B_1$ and $B_2$ be $\gl_{\infty}$-crystals. A {\it morphism}
$\psi : B_1 \rightarrow B_2$ is a map from $B_1\cup\{{\bf 0}\}$ to
$B_2\cup\{{\bf 0}\}$ such that
\begin{itemize}
\item[(1)] $\psi(\bf{0})=\bf{0}$,

\item[(2)] ${\rm wt}(\psi(b))={\rm wt}(b)$,
$\varepsilon_i(\psi(b))=\varepsilon_i(b)$, and
$\varphi_i(\psi(b))=\varphi_i(b)$ whenever $\psi(b)\neq \bf{0}$,

\item[(3)] $\psi(\te_i b)=\te_i\psi(b)$ for $b\in B_1$ such that $\psi(b)\neq \bf{0}$ and
$\psi(\te_i b)\neq \bf{0}$,

\item[(4)] $\psi(\tf_i
b)=\tf_i\psi(b)$ for $b\in B_1$ such that $\psi(b)\neq \bf{0}$ and
$\psi(\tf_i b)\neq \bf{0}$.
\end{itemize}
We call $\psi$ an {\it embedding} and $B_1$ a {\it subcrystal of}
$B_2$ when $\psi$ is injective, and {\it strict} if $\psi :
B_1\cup\{{\bf 0}\} \rightarrow B_2\cup\{{\bf 0}\}$ commutes with
$\te_i$ and $\tf_i$ ($i\in\Z$), where we assume that $\te_i{\bf 0}=\tf_i{\bf 0}={\bf 0}$.

We define the {\it tensor product  of} $B_1$ and $B_2$ to be the set
$B_1\otimes B_2=\{\,b_1\otimes b_2\,|\,b_i\in B_i\,\, (i=1,2)\,\}$
with {\allowdisplaybreaks
\begin{equation*}
{\rm wt}(b_1\otimes b_2)={\rm wt}(b_1)+{\rm wt}(b_2),
\end{equation*}
\begin{equation*}
\varepsilon_i(b_1\otimes b_2)= {\rm
max}(\varepsilon_i(b_1),\varepsilon_i(b_2)-\langle {\rm
wt}(b_1),h_i\rangle),
\end{equation*}
\begin{equation*}
\varphi_i(b_1\otimes b_2)= {\rm max}(\varphi_i(b_1)+\langle {\rm
wt}(b_2),h_i\rangle,\varphi_i(b_2)),
\end{equation*}
\begin{equation*}
{\te}_i(b_1\otimes b_2)=
\begin{cases}
{\te}_i b_1 \otimes b_2, & \text{if $\varphi_i(b_1)\geq \varepsilon_i(b_2)$}, \\
b_1\otimes {\te}_i b_2, & \text{if
$\varphi_i(b_1)<\varepsilon_i(b_2)$},
\end{cases}
\end{equation*}
\begin{equation*}
{\tf}_i(b_1\otimes b_2)=
\begin{cases}
{\tf}_i b_1 \otimes b_2, & \text{if  $\varphi_i(b_1)>\varepsilon_i(b_2)$}, \\
b_1\otimes {\tf}_i b_2, & \text{if $\varphi_i(b_1)\leq
\varepsilon_i(b_2)$},
\end{cases}
\end{equation*}}
\noindent for $i\in\Z$, where we assume that ${\bf 0}\otimes
b_2=b_1\otimes {\bf 0}={\bf 0}$.

 For $b_i\in B_i$ ($i=1,2$), let $C(b_i)$ denote the
connected component of $b_i$ in $B_i$ as a $\Z$-colored oriented
graph. We say that $b_1$ {\it is equivalent to} $b_2$ if there is an
isomorphism of crystal graphs $C(b_1)\rightarrow C(b_2)$ sending
$b_1$ to $b_2$.

Let $\B$ be a $\gl_{\infty}$-crystal given by
\begin{equation*}
\cdots \ \stackrel{-2}{\longrightarrow} -2 \
\stackrel{-1}{\longrightarrow} -1 \stackrel{0} {\longrightarrow} 1
\stackrel{1}{\longrightarrow} 2\stackrel{2}{\longrightarrow} \cdots,
\end{equation*}
where ${\rm wt}(k)=\epsilon_k$, and $\varepsilon_{i}(k)$ (resp.
$\varphi_{i}(k)$) is the number of $i$-colored arrows coming into
$b$ (resp. going out of $k$) for $k\in\B$.

Let $\gl_{< 0}$ and $\gl_{>0}$ be the subalgebras of $\gl_{\infty}$
spanned by $\{\,E_{ij}\,|\,i,j\in\Z_{< 0}\,\}$ and
$\{\,E_{ij}\,|\,i,j\in\Z_{> 0}\,\}$, respectively. We can define
$\gl_{<0}$-crystals (resp. $\gl_{>0}$-crystals) as in Definition
\ref{crystal graph} with respect to $\te_i$'s and $\tf_i$'s for
$i\in\Z_{<0}$ (resp. $i\in\Z_{>0}$), and view
$\B_{>0}=\{\,k\in\B\,|\,k>0\,\}$ as a
 $\gl_{>0}$-crystal  and $\B_{<0}=\{\,k\in\B\,|\,k<0\,\}$ as a $\gl_{<0}$-crystal.
 In addition, let us consider a $\gl_{<0}$-crystal $\B_{<0}^\vee$
given as follows;
\begin{equation*}
-1^\vee \stackrel{-1}{\longrightarrow}
-2^\vee\stackrel{-2}{\longrightarrow}-3^\vee\stackrel{-3}{\longrightarrow}
\cdots,
\end{equation*}
where ${\rm wt}(-k^\vee)=-\epsilon_{-k}$ for $k>0$. Note that
$\B_{<0}^\vee$ is the dual crystal of $\B_{<0}$ (cf.\cite{Kas94}).

For $\lambda\in P$, let $T_{\lambda}=\{\,t_{\lambda}\,\}$ be a
$\gl_{\infty}$-crystal with ${\rm wt}(t_{\lambda})=\lambda$, $\te_i
t_{\lambda}=\tf_{i} t_{\lambda}={\bf 0}$, and
$\varepsilon_{i}(t_{\lambda})=\varphi_{i}(t_{\lambda})=-\infty$ for
$i\in\Z$.

\subsection{Semistandard tableaux}
Let $\mathscr{P}$ be the set of partitions. We identify a partition
$\lambda$ with a {\it Young diagram} or a subset
$\{\,(i,j)\,|\,1\leq j\leq \lambda_i\,\}$ of
$\mathbb{N}\times\mathbb{N}$ (cf. \cite{Mac95}). The number of
non-zero parts of $\lambda$ is denoted by $\ell(\lambda)$, called
the {\it length of $\lambda$}. For $\lambda=(\lambda_i)_{i\geq
1}\in\cP$, $\lambda'=(\lambda'_i)_{i\geq 1}$ denotes the {\it
conjugate of $\lambda$}. For $\mu\in\cP$ with $\mu\subset \lambda$,
$\lambda/\mu$ denotes the {\it skew Young diagram}, and
$|\lambda/\mu|$ denotes the number of boxes in the diagram. We
denote by $\lambda^\pi$ the skew Young diagram obtained by
$180^{\circ}$-rotation of $\lambda$, which is called of {\it
anti-normal} shape. Put
$\mathscr{P}^{\pi}=\{\,\lambda^\pi\,|\,\lambda\in\cP\,\}$, the set
of anti-normal shaped skew Young diagrams.

Let $\A$ be a linearly ordered set. For a skew Young diagram
$\lambda/\mu$, a tableau $T$ obtained by filling $\lambda/\mu$ with
entries in $\A$ is called a {\it semistandard tableau of shape
$\lambda/\mu$} if the entries in each row are weakly increasing from
left to right, and the entries in each column are strictly
increasing from top to bottom. We write ${\rm sh}(T)=\lambda/\mu$.
We denote by $SST_\A(\lambda/\mu)$ the set of all semistandard
tableaux of shape $\lambda/\mu$ with entries in $\A$.
Let $\W_{\A}$ be the set of finite words in $\A$. We associate to each $T\in
SST_\A(\lambda/\mu)$ a word $w(T)\in\W_\A$ which is obtained by
reading the entries of $T$ row by row from top to bottom, and from
right to left in each row.

We suppose that $\A=\B$, $\B_{> 0}$, $\B_{< 0}$, and $\B^\vee_{<0}$,
with the linear ordering $<$ induced from the partial ordering on
$P$.  Then
$\W_{\B}$ is a $\gl_{\infty}$-crystal since we may view each
non-empty finite word $w=w_1\cdots w_r$ as $w_1\otimes\cdots\otimes
w_r\in \B^{\otimes r}$. Similarly, $\W_{\B_{>0}}$ (resp.
$\W_{\B_{<0}}$ or $\W_{\B_{<0}^\vee}$) becomes a $\gl_{>0}$-crystal
(resp. $\gl_{<0}$-crystal).
Sending $T$ to $w(T)$ gives an injective
map from $SST_\A(\lambda/\mu)$ to $\W_{\A}$, and the image of
$SST_\A(\lambda/\mu)$ together with $\{{\bf 0}\}$ is invariant under
the operators $\te_i,\tf_i$. Hence it is a
crystal graph \cite{KN}.

Suppose that $\lambda/\mu\in \mathscr{P}$ or $\mathscr{P}^\pi$. Then
$SST_\A(\lambda/\mu)$ is connected. In particular, if $\A=\B_{>0}$
(resp. $\B_{<0}^\vee$), then $SST_\A(\lambda /\mu)$ contains a
highest weight element $H_{\lambda/\mu}$, where in each $k$th column
of $H_{\lambda/\mu}$, the $l$th entry from the top position is
filled with $l$ (resp. $-l^\vee$).

\section{The RSK correspondence and $\gl_{\infty}$-crystals}\label{RSK}
In this section, we define two $\gl_{\infty}$-crystal structures on
the set $\M$ of matrices with non-negative integral entries of
finite support and the set $\N$ of pairs of semistandard tableaux of
the same shape. We show that the  RSK correspondence, which is
a bijection from $\M$ to $\N$, is an isomorphism of
$\gl_{\infty}$-crystals.

Since it is already known that the RSK correspondence is a morphism
of $(\gl_{<0},\gl_{>0})$-bicrystals \cite{DK,La}, our main result in
this section is to extend it as a $\gl_{\infty}$-crystal morphism by defining the missing operators
$\te_0,\tf_0$ on $\M$ and $\N$, which are compatible with the RSK
algorithm.

\subsection{Crystal of integral matrices}
Let
\begin{equation}
\begin{split}
\Omega=\{\,& (\bi,\bj)\in
\W_{\B_{<0}^\vee}\times \W_{\B_{>0}}\,| \\
& \text{(1) $\bi=i_1\cdots i_r$ and $\bj=j_1\cdots j_r$  for some $r\geq 0$}, \\
& \text{(2) $(i_1,j_1)\leq \cdots \leq (i_r,j_r)$} \},
\end{split}
\end{equation}
where for $(i,j)$ and $(k,l)\in \B_{<0}^\vee\times \B_{>0}$,
\begin{equation*}\label{partial order}
(i,j)< (k,l) \ \ \ \ \Longleftrightarrow \ \ \ \
\begin{cases}
j<l & \text{or}, \\
j=l\ \text{and} \ i>k.
\end{cases}
\end{equation*}
Similarly, let $\Omega'$ be the set of pairs $(\bk,\bl)\in
\W_{\B_{<0}^\vee}\times \W_{\B_{>0}}$ such that (1) $\bk=k_1\cdots
k_r$ and $\bl=l_1\cdots l_r$  for some $r\geq 0$ and (2)
$(k_1,l_1)\leq' \cdots \leq' (k_r,l_r)$, where
\begin{equation*}\label{partial order}
(i,j)<' (k,l) \ \ \ \ \Longleftrightarrow \ \ \ \
\begin{cases}
i<k & \text{or}, \\
i=k\ \text{and} \ j>l.
\end{cases}
\end{equation*}
Then $\Omega$ is a $\gl_{<0}$-crystal, where $x_i
(\bi,\bj)=(x_i\bi,\bj)$ for $(\bi,\bj)\in \Omega$, $x=\te,\tf$ and
$i\in \Z_{<0}$. Here, we assume that $x_i (\bi,\bj)={\bf 0}$ if
$x_i\bi={\bf 0}$. Similarly, $\Omega'$ is a $\gl_{>0}$-crystal,
where $x_j(\bk,\bl)=(\bk,x_j\bl)$ for $(\bk,\bl)\in \Omega'$,
$x=\te,\tf$ and $j\in \Z_{>0}$.

Consider the following set of $\mathbb{N}\times\mathbb{N}$ matrices
with non-negative integers of finite support;
\begin{equation}
\begin{split}
\M=& \{\,A=(a_{-i^\vee, j})_{i,j\geq 1}\,|\, a_{-i^\vee,
j}\in\mathbb{Z}_{\geq 0},\ \sum_{i,j\geq 1}a_{-i^\vee, j}<\infty\,
\}.
\end{split}
\end{equation}

For $(\bi,\bj)\in \Omega$, define $A(\bi,\bj)=(a_{-i^\vee,j})$ to be
the matrix in $\M$, where $a_{-i^\vee,j}$ is the number of $k$'s
such that $(i_k,j_k)=(-i^\vee,j)$. Then the map $(\bi,\bj)\mapsto
A(\bi,\bj)$ is a bijection between $\Omega$ and $\M$, where the pair
of empty words $(\emptyset,\emptyset)$ corresponds to zero matrix,
say $\mathbb{O}$. Similarly, we have a bijection $(\bk,\bl)\mapsto
A(\bk,\bl)$ from $\Omega'$ to $\M$.

With these bijections, $\M$ becomes a crystal graph for both
$\gl_{<0}$ and $\gl_{>0}$.
Moreover, the operators $\te_i,\tf_i$ commute with $\te_j,\tf_j$ for
$i\in\mathbb{Z}_{<0}$ and $j\in\mathbb{Z}_{>0}$, and hence $\M$
becomes a $(\gl_{<0},\gl_{>0})$-bicrystal (cf.\cite{DK,La}).

Now, for $A=(a_{-i^\vee, j})\in\M$, we define
\begin{equation}
\begin{split}
&\te_0 A=
\begin{cases}
A-E_{-1^\vee,1}, & \text{if $a_{-1^\vee,1}\neq 0$}, \\
{\bf 0}, & \text{otherwise},
\end{cases}\\
&\tf_0 A= A+E_{-1^\vee,1},
\end{split}
\end{equation}
where $E_{-1^\vee,1}\in \M$ denotes the elementary matrix with 1 at
the position  $(-1^\vee,1)$ and $0$ elsewhere. Put ${\rm
wt}(A)=\sum_{i,j>0} a_{-i^\vee\, j}(-\epsilon_{-i}+\epsilon_{j})$,
$\varepsilon_0(A)=a_{-1^\vee,1}$, and $\varphi_0(A)=\langle {\rm
wt}(A),h_0 \rangle + \varepsilon_0(A)$. Then we have the following.

\begin{prop}\label{crystalM}
$\M$ is a $\gl_{\infty}$-crystal, and
\begin{equation*}
\M=\{\,\tf_{i_1}\cdots\tf_{i_r}{\mathbb{O}}\,|\,r\geq 0, \
i_1,\ldots,i_r\in \Z\,\}\setminus\{{\bf 0}\}.
\end{equation*}
In particular, $\M$ is connected with highest weight element
$\mathbb{O}$.
\end{prop}
\pf It is easy to see that $\M$ is a $\gl_{\infty}$-crystal. Let
$A\in \M$ be given. We claim that
$A=\tf_{i_1}\cdots\tf_{i_r}{\mathbb{O}}$ for some $r\geq 0$ and
$i_1,\ldots,i_r\in\Z$. We use induction on
$s(A)=\sum_{i,j}a_{-i^\vee,j}$. If $s(A)=0$, then it is clear. Let
$s(A)$ be positive. First, $A$ is connected to a diagonal matrix
$A^\circ=(a^\circ_{-i^\vee,j})$ such that $a^\circ_{-1^\vee,1}\geq
a^\circ_{-2^\vee,2}\geq a^\circ_{-3^\vee,3}\geq \ldots$ since $\M$
is a $(\gl_{<0},\gl_{>0})$-bicrystal (cf.\cite{DK,K07}). That is,
$\te_{j_1}\cdots\te_{j_r}A=A^\circ$ for some $r\geq 0$ and
$j_1,\ldots,j_r\in\Z^{\times}$ and $\te_iA^\circ={\bf 0}$ for all
$i\in\Z^{\times}$. If $a^\circ_{-1^\vee,1}=0$, then
$A^\circ=\mathbb{O}$. If not, then $\te_0 A^\circ\neq {\bf 0}$ and
$s(A^\circ)=s(A)-1$. Hence, the proof completes by induction
hypothesis. \qed

\subsection{Crystal of bitableaux}
By the RSK algorithm, each $A\in\M$ is in one-to-one correspondence
with $(P(A),Q(A))$ in $SST_{\B^\vee_{<0}}(\lambda)\times
SST_{\B_{>0}}(\lambda)$ for some $\lambda\in\mathscr{P}$ \cite{Kn}.
In what follows, we need a variation of this correspondence with
anti-normal shaped tableaux. Let us describe it in detail.

Let $\nu\in\mathscr{P}$ and $T\in SST_\A(\nu^\pi)$ be given. For
$a\in \A$, we define $T \leftarrow a$ to be the tableau of an
anti-normal shape obtained from $T$ by applying the following
procedure; (1) let $a'$ be the largest entry in the right-most
column which is smaller than or equal to $a$, (2) replace {$a'$} by
{$a$}. If there is no such $a'$, put {$a$} at the top of the column
and stop the procedure, (3) repeat (1) and (2) on the next column
with {$a'$}. For $w=w_1\ldots w_r\in\W_\A$, we define ${\bf P}(w)$
to be $(\cdots(({w_1}\leftarrow{w_2}\,)\leftarrow{w_3}\,)
\cdots)\leftarrow{w_r}$. Note that $w$ is equivalent to ${\bf P}(w)$
as elements of crystals.

Let $A\in \M$ be given with $A=A(\bi,\bj)=A(\bk,\bl)$ for
$(\bi,\bj)\in \Omega$ and $(\bk,\bl)\in\Omega'$. Let $\bi_{\rm rev}$
and $\bl_{\rm rev}$ be the reverse word of $\bi$ and $\bl$,
respectively. We define
\begin{equation}
{\bf P}(A)={\bf P}(\bi_{\rm rev}), \ \ \ {\bf Q}(A)={\bf P}(\bl_{\rm
rev}).
\end{equation}

Let $\bi=i_1\ldots i_r$. For $1\le k \le r$, let us fill the box
with $c$ if it is created when $i_k$ is inserted into
$(({i_r}\leftarrow{i_{r-1}}\,) \cdots)\leftarrow{i_{k+1}}$ and
$j_k=c$. Then we have a tableau $Q(\bi_{\rm rev})\in
SST_{\B_{>0}}(\nu^\pi)$ with $\nu^\pi={\rm sh}P(\bi_{\rm
rev})\in\mathscr{P}^\pi$. By the symmetry of the RSK correspondence,
we have $Q(\bi_{\rm rev})={\bf P}(\bl_{\rm rev})={\bf Q}(A)$. Put
\begin{equation}
\N = \bigsqcup_{\nu \in\mathscr{P}}SST_{\B^\vee_{<0}}(\nu^\pi)\times
SST_{\B_{>0}}(\nu^\pi).
\end{equation}
Hence we have a bijection
\begin{equation}\label{kappa}
\kappa : \M \longrightarrow \N,
\end{equation}
where $\kappa(A)=({\bf P}(A),{\bf Q}(A))$.

\begin{ex}\label{RSK ex}{\rm
Let $$A=(a_{-i^\vee,j})_{1\leq i,j\leq 3}=\left(
         \begin{array}{ccc}
           1 & 0 & 1 \\
           2 & 1 & 0 \\
           0 & 2 & 0 \\
         \end{array}
       \right)
\in\M,$$ where we assume that $a_{-i^\vee, j}=0$ unless $1\leq
i,j\leq 3$. Then
$$\bi=-2^\vee \!\!-2^\vee\!\! -1^\vee\!\! -3^\vee\!\! -3^\vee\!\! -2^\vee\!\! -1^\vee,
\ \ \ \bl=3\ 1\ 2\ 1\ 1\ 2\ 2,$$ and
$$
{\bf P}(A)=
  \begin{array}{cccc}
    \!\!\!\! & -1^\vee\!\!\!\! & -2^\vee\!\!\!\! & -2^\vee \\
    -1^\vee\!\!\!\! & -2^\vee\!\!\!\! & -3^\vee\!\!\!\! & -3^\vee \\
  \end{array} \ \ , \ \
{\bf Q}(A)=
  \begin{array}{cccc}
     & 1 \ & 1\ & 1\ \\
    2\ & 2\ & 2\ & 3\ \\
  \end{array}\ \ .
$$
 }
\end{ex}

Clearly $\N$ is a $(\gl_{<0},\gl_{>0})$-bicrystal, and for $A\in
\M$, ${\bf P}(A)$ (resp. ${\bf Q}(A)$) is equivalent to $A$ as
elements of $\gl_{<0}$ (resp. $\gl_{>0}$)-crystals. Summarizing, we
have the following.

\begin{prop}\label{bicrystal iso}
$\kappa$ is a $(\gl_{<0},\gl_{>0})$-bicrystal isomorphism.\qed
\end{prop}

Now, let us describe the $\gl_{\infty}$-crystal structure on $\N$.
Suppose that $(S,T)\in\N$ is given. For each $k$th column of $\nu$
enumerated from right to left, let $s_k$ and $t_k$ be the smallest
entries in the $k$th column of $S$ and $T$, respectively. We assign
\begin{equation}\label{signs}
\sigma_k=
\begin{cases}
+ \ , & \text{if the $k$th column is empty}, \\
+ \ ,& \text{if $s_k>-1^\vee$ and $t_k>1$}, \\
- \ ,& \text{if $s_k=-1^\vee$ and $t_k=1$,}\\
\ \cdot \ \, ,& \text{otherwise}.
\end{cases}
\end{equation}

In the sequence $\sigma=(\ldots\sigma_2,\sigma_1)$, we replace a
pair $(\sigma_{s'},\sigma_{s})=(-,+)$, where $s'>s$ and
$\sigma_t=\cdot$ for $s<t<s'$, by $(\,\cdot\,,\,\cdot\,)$, and repeat
this process as far as possible until we get a sequence with no $-$
placed to the left of $+$.  We call this sequence the {\it
$0$-signature} of $(S,T)$

We call the left-most $-$ in the $0$-signature of $(S,T)$ the {\it
$0$-good $-$ sign}, and define $\te_0 (S,T)$ to be the bitableaux
obtained from $(S,T)$  by removing $-1^\vee$ and $1$ in the corresponding columns.
If there is no $0$-good $-$ sign, then we define $\te_0 (S,T)={\bf
0}$. We call the right-most $+$ in the $0$-signature of $(S,T)$ the
{\it $0$-good $+$ sign}, and define $\tf_0 (S,T)$ to be the
bitableaux obtained from $(S,T)$ by adding $-1^\vee$ and $1$ on top of the
corresponding columns. If there is no $0$-good $+$ sign, then we
define $\tf_0 (S,T)={\bf 0}$.

\begin{ex}\label{Nex}{\rm Let $(S,T)\in\N$ be given with  $\sigma$ as follows.
\begin{equation*}
\begin{split}
S=
\begin{array}{cccc}
   &  & -1^\vee \!\!\! & -3^\vee \!\!\! \\
  -1^\vee \!\!\! & -2^\vee \!\!\! & -4^\vee \!\!\! & -4^\vee \!\!\!
\end{array}\ \ , \ \
T=
\begin{array}{cccc}
     &  & 1\ & 3\ \\
  1\ & 1\ & 3\ & 4\
\end{array}\ \ ,\\
\sigma=(\ldots\sigma_5,\sigma_4,\sigma_3,\sigma_2,\sigma_1)=(\ldots
+, +, -,\, \cdot \,, -, +).
\end{split}
\end{equation*}
Then the $0$-good $-$ sign is $\sigma_4$, and $0$-good $+$ sign is
$\sigma_5$. Hence,
\begin{equation*}
\begin{split}
\te_0(S,T)&=\left(
\begin{array}{ccc}
     & -1^\vee \!\!\! & -3^\vee \!\!\! \\
     -2^\vee \!\!\! & -4^\vee \!\!\! & -4^\vee \!\!\!
\end{array}\  ,  \
\begin{array}{ccc}
       & 1\ & 3\ \\
     1\ & 3\ & 4\
\end{array} \right) \ \ ,\\
\tf_0(S,T)&=\left(
\begin{array}{ccccc}
   & & & -1^\vee\!\!\! & -3^\vee\!\!\! \\
 -1^\vee\!\!\! & -1^\vee\!\!\! & -2^\vee\!\!\! & -4^\vee\!\!\! & -4^\vee\!\!\!
\end{array}\  ,  \
\begin{array}{ccccc}
    &  &  & 1\ & 3\ \\
  1\ & 1\ & 1\ & 3\ & 4\
\end{array} \right) \ \ .
\end{split}
\end{equation*} }
\end{ex}

\begin{prop}
$\N$ is a $\gl_{\infty}$-crystal, and
\begin{equation*}
\N=\{\,\tf_{i_1}\cdots\tf_{i_r}{(\emptyset,\emptyset)}\,|\,r\geq 0,
\ i_1,\ldots,i_r\in \Z\,\}\setminus\{{\bf 0}\}.
\end{equation*}
In particular, $\N$ is connected with highest weight element
$(\emptyset,\emptyset)$.
\end{prop}
\pf It is straightforward to check that $\te_0(S,T), \tf_0(S,T) \in
\N\cup \{{\bf 0}\}$. Let $(S,T)$ be given with ${\rm sh}(S)={\rm
sh}(T)$ non-empty. We may assume that $\te_i(S,T)={\bf 0}$ for all
$i\in \Z^{\times}$. Then $S$ (resp. $T$) is a highest weight element
of a $\gl_{<0}$-crystal  (resp. $\gl_{>0}$-crystal), where in each
$k$th column of $S$ (resp. $T$), the $l$th entry from the top
position is filled with $-l^\vee$ (resp. $l$). Hence $\te_0(S,T)\neq
{\bf 0}$. If we use induction on $|{\rm sh}(S)|$, then we conclude
that $\tf_{i_1}\cdots\tf_{i_r}{(\emptyset,\emptyset)}=(S,T)$ for
some $r\geq 0$ and $i_1,\ldots,i_r\in \Z$. \qed

\subsection{Isomorphism}
Now we are in a position to state the main result in this section.

\begin{thm}\label{RSK iso}
The map $\kappa :\M \longrightarrow \N$ is a $\gl_{\infty}$-crystal isomorphism.
\end{thm}
\pf By Proposition \ref{bicrystal iso}, it suffices to show that
$\kappa$ commutes with $\te_0$ and $\tf_0$. More precisely, we claim
that for $A\in\M$ and $k\geq 1$,
\begin{equation}\label{mainstep}
\kappa(\tf_0^kA)=\kappa(kE_{-1^\vee,1}+A)=\tf_0^k\kappa(A).
\end{equation}
We use induction on $t(A)=\sum_{k\geq 1}a_{-k^\vee,1}$. We may
assume that $a_{-1^\vee,1}=0$.

If $t(A)=0$, then it is not difficult to see that \eqref{mainstep}
holds. We suppose that $t(A)>0$. Let $-p^\vee$ be the largest one
such that $a_{-p^\vee,1}\neq 0$. Let $B=A-E_{-p^\vee,1}$. By
induction hypothesis, we have for $k\geq 1$
$$\kappa(\tf_0^kB)=\kappa(kE_{-1^\vee,1}+B)=\tf_0^k\kappa(B).$$  For $k\geq 1$, let $\kappa(kE_{-1^\vee,1}+X)=({\bf
P}^{(k)}(X),{\bf Q}^{(k)}(X))$  with $X=A,B$. Then ${\bf
P}^{(k)}(A)={\bf P}^{(k)}(B)\leftarrow -p^\vee$ by definition of
${\bf P}(\cdot)$ and ${\bf Q}^{(k)}(A)$ is obtained from ${\bf
Q}^{(k)}(B)$ by filling the corresponding box, say $c$, in ${\rm
sh}({\bf P}^{(k)}(A))/{\rm sh}({\bf P}^{(k)}(B))$ with $1$. For
convenience, let us write
$\kappa(kE_{-1^\vee,1}+A)=\kappa(kE_{-1^\vee,1}+B)\leftarrow
(-p^\vee,1)$. Let $$\sigma=(\ldots,\sigma_2,\sigma_1),\ \ \
\sigma'=(\ldots,\sigma'_2,\sigma'_1)$$ be the sequences of signs
associated with $\kappa(kE_{-1^\vee,1}+B)$ and
$\kappa(kE_{-1^\vee,1}+A)$, respectively (see (\ref{signs})), and
let
$$\tilde{\sigma}=(\ldots,\tilde{\sigma}_2,\tilde{\sigma}_1), \ \ \
\widetilde{\sigma'}=(\ldots,\widetilde{\sigma'}_2,\widetilde{\sigma'}_1)$$
be the $0$-signatures of $\kappa(kE_{-1^\vee,1}+B)$ and
$\kappa(kE_{-1^\vee,1}+A)$, respectively.

Suppose that by the insertion of $-p^\vee$ into ${\bf P}^{(k)}(B)$,
$c$ is filled with $-q^\vee$ for some $q\leq p$, and it is located
at the $t$th column  enumerated from the rightmost one.

\textsc {Case 1}. $q>1$. Let $kE_{-1^\vee,1}+B=A(\bi,\bj)$ with
$(\bi,\bj)\in\Omega$. Consider the horizontal strip made by
inserting the subwords of $\bi_{\rm rev}$ corresponding to the first
column of $kE_{-1^\vee,1}+B$.

Then we observe the following facts;
\begin{itemize}
\item[(1)] no $-1^\vee$ has been bumped out in the
bumping path for ${\bf P}^{(k)}(B)\leftarrow -p^\vee$.

\item[(2)] by
induction hypothesis  all $-1^\vee$'s which have been added on ${\bf
P}(B)$ by applying $\tf^k_0$ to $\kappa(B)$ are placed to the right
of $-q^\vee$ in the $t$th column, and they do not intersect with the
bumping path for ${\bf P}^{(k)}(B)\leftarrow -p^\vee$.

\item[(3)] the insertion of $-p^\vee$
into ${\bf P}^{(k)}(B)$ does not change the sign $\sigma_k$  for
$1\leq k\leq t-1$, and hence $\sigma_k=\sigma'_k$ for $1\leq k\leq
t-1$.

\end{itemize}

Hence we have $$\te_0^k[\kappa(kE_{-1^\vee,1}+B)\leftarrow
(-p^\vee,1)]=\kappa(B)\leftarrow (-p^\vee,1)=\kappa(A),$$ and
\begin{equation*}
\tf_0^k\kappa(A)=\kappa(kE_{-1^\vee,1}+B)\leftarrow
(-p^\vee,1)=\kappa(kE_{-1^\vee,1}+A).
\end{equation*}

\textsc{ Case 2}. $q=1$. Consider the bumping path for ${\bf
P}^{(k)}(B)\leftarrow -p^\vee$. Then there exists $1\leq s\leq t$
such that
\begin{itemize}
\item[(1)] $-x^\vee$ ($x\geq 2$) has been bumped out from the $(k-1)$th column
and placed at the $k$th column for $2\leq k\leq s$,

\item[(2)] $-1^\vee$ has
been bumped out from the $(k-1)$th column and placed at the $k$th
column for $s+1\leq k\leq t$.
\end{itemize}

As in \textsc{Case 1}, it follows that all $-1^\vee$'s which have
been added to ${\bf P}(B)$ by applying $\tf^k_0$ to $\kappa(B)$ are
placed to the right of the $t$th column, and $\sigma_r=\sigma'_r$
for $1\leq r\leq s$.

Since all  $-1^\vee$'s in the $r$th column of ${\bf P}^{(k)}(B)$ for
$s\leq r \leq t-1$ have been shifted to the left by one column by
the insertion of $-p^\vee$ to ${\bf P}^{(k)}(B)$, we have
$\sigma_r=\sigma'_r$ for $s+1\leq r\leq t-1$. Note that
${\sigma}_t=+$ and $\sigma'_t=-$.

Let $u$ be the top entry of the $s$th column in ${\bf Q}^{(k)}(B)$.
If $u=1$, then we have $\sigma_s=-$ and $\sigma'_s=\cdot$. If $u>1$,
then we have $\sigma_s=\cdot$ and $\sigma'_s=+$. Now, comparing
$\sigma$ and $\sigma'$ (hence $\tilde{\sigma}$ and
$\tilde{\sigma}$), it is not difficult to see that
\begin{equation*}
\te_0^k\kappa(kE_{-1^\vee,1}+A)=\te_0^k[\kappa(kE_{-1^\vee,1}+B)\leftarrow
(-p^\vee,1)]=\kappa(B)\leftarrow (-p^\vee,1)=\kappa(A).
\end{equation*}
This completes the proof. \qed

\begin{ex}{\rm Consider
$$\kappa
\left(
  \begin{array}{cccc}
    1 & 0 & 1 & 0 \\
    0 & 0 & 0 & 1 \\
    1 & 0 & 0 & 0 \\
    1 & 0 & 1 & 0 \\
  \end{array}
\right)=\left(
\begin{array}{cccc}
   &  & -1^\vee\!\!\! & -3^\vee\!\!\! \\
  -1^\vee\!\!\! & -2^\vee\!\!\! & -4^\vee\!\!\! & -4^\vee\!\!\!
\end{array}\ \ , \ \
\begin{array}{cccc}
     &  & 1\ & 3\ \\
  1\ & 1\ & 3\ & 4\
\end{array}\right)\ \ .
$$
Applying $\te_0$ and $\tf_0$ on both sides, we get
$$
\kappa \left(
  \begin{array}{cccc}
    {\bf 0} & 0 & 1 & 0 \\
    0 & 0 & 0 & 1 \\
    1 & 0 & 0 & 0 \\
    1 & 0 & 1 & 0 \\
  \end{array}
\right)=\left(
\begin{array}{cccc}
   &  & -1^\vee\!\!\! & -3^\vee\!\!\! \\
    & -2^\vee\!\!\! & -4^\vee\!\!\! & -4^\vee\!\!\!
\end{array}\ \ , \ \
\begin{array}{cccc}
     &  & 1\ & 3\ \\
    & 1\ & 3\ & 4\
\end{array}\right)\ \ ,
$$ and
$$
\kappa \left(
  \begin{array}{cccc}
   {\bf 2} & 0 & 1 & 0 \\
    0 & 0 & 0 & 1 \\
    1 & 0 & 0 & 0 \\
    1 & 0 & 1 & 0 \\
  \end{array}
\right)=\left(
\begin{array}{ccccc}
  & &  & -1^\vee\!\!\! & -3^\vee\!\!\! \\
 -1^\vee\!\!\! & -1^\vee\!\!\! & -2^\vee\!\!\! & -4^\vee\!\!\! & -4^\vee\!\!\!
\end{array}\ \ , \ \
\begin{array}{ccccc}
   &  &  & 1\ & 3\ \\
 1\ & 1\ & 1\ & 3\ & 4\
\end{array}\right)\ \ ,
$$
respectively (see Example \ref{Nex}). }
\end{ex}

\section{Crystal graphs of generalized Verma modules}
Let $\mathfrak{u}_{\pm}$ be the subalgebra of $\gl_{\infty}$ spanned
by $E_{ij}$ for $i<0$, $j>0$ (resp. $i>0$, $j<0$). Let
$\frak{p}=\gl_{<0}\oplus\gl_{>0}\oplus \frak{u}_+$ be a maximal
parabolic subalgebra. Then we have
$\gl_{\infty}=\frak{u}_-\oplus\frak{p}$. The set of roots for the
nilradical $\mathfrak{u}_-$ is given by
$\Delta(\mathfrak{u}_-)=\{\,-\epsilon_i+\epsilon_{j}\,|\,i>0,j<0\,\}$.
Let $U(\mathfrak{u}_-)$ be the universal enveloping algebra of
$\mathfrak{u}_-$. By PBW theorem, $U(\mathfrak{u}_-)$ has a basis
parameterized by $\M$.

In this section, we prove that the $\gl_{\infty}$-crystal $\M$ is
the crystal graph of a generalized Verma module $U(\mathfrak{u}_-)$
or its $q$-analogue (cf.\cite{Br}) in the sense that it is the limit
of the crystal graphs of the integrable highest weight
$\gl_{\infty}$-modules with highest weight $n\Lambda_0$ as
$n\rightarrow \infty$.

\subsection{Crystal $\B(n\Lambda_0)$ }
Let $\F$ be the set of {\it semi-infinite words} $$w=\cdots
w_{-3}w_{-2}w_{-1}$$ with letters in $\B$ such that
\begin{itemize}
\item[(1)] $w_{i-1}< w_i$ for all $i< 0$,

\item[(2)] there exists an integer $c\in\mathbb{Z}$ such that $w_i=i+c$ for $i\ll
0$.
\end{itemize}
For $w\in \F$, we define ${\rm
wt}(w)=\Lambda_0+\sum_{k\in\B}m_k\epsilon_k\in {P}$, where
$m_k=\bigl|\{\,i\,|\,w_i=k\ \,\}\bigr|-\delta_{-k,|k|}$. It is
well-defined since $m_k=0$ for almost all $k\in \B$. For each $i\in
\Z$, we define the operators $\te_i,\tf_i : \F \longrightarrow
\F\cup\{{\bf 0}\}$ by the same way as we do on $\W_{\B}$. Then
$\te_i$ and $\tf_i$ are well-defined, and $\F$ is a
$\gl_{\infty}$-crystal. For $i\leq 0$, let $H_{\Lambda_i}= \ldots \
i-3\ i-2 \ i-1 $, and for $i>0$, let $H_{\Lambda_i}= \ldots  -2
-1\,1 \ldots\, i-1 \ i$. We have the following decomposition as
$\gl_{\infty}$-crystals
$$\F=\bigsqcup_{i\in\mathbb{Z}}\B(\Lambda_i),$$ where $\B(\Lambda_i)$
is the connected component of $H_{\Lambda_i}$ with ${\rm
wt}(H_{\Lambda_i})=\Lambda_i$. Recall that $\F$ is the crystal graph
of the Fock space representation, which can be realized as the space
of semi-infinite wedge vectors, and $\B(\Lambda_i)$ is the crystal
graph of the irreducible highest weight $\gl_{\infty}$-module with
highest weight $\Lambda_i$ (cf.\cite{MM}).

Let $\lambda=(\lambda_1,\cdots,\lambda_n)$ be a sequence of
non-increasing $n$ integers, called a {\it generalized partition of
length $n$}. We call an $n$-tuple of semi-infinite words
$\w=(w^{(1)},\cdots,w^{(n)})$ a {\it semi-infinite semistandard
tableau of shape $\lambda$} if
\begin{itemize}
\item[(1)] $w^{(i)}=\cdots w^{(i)}_{-3}w^{(i)}_{-2}w^{(i)}_{-1}\in\B(\Lambda_{\lambda_{n-i+1}})$ for $1\leq i\leq n$,

\item[(2)] $w^{(i+1)}_{k+d_i}\leq w^{(i)}_k$ for $1\leq i<n$ and $k< 0$,
where $d_i=\lambda_{n-i+1}-\lambda_{n-i}$.
\end{itemize}

We may identify each $\w$ with a semistandard tableau with
infinitely many rows and $n$ columns, where each row of $\w$ reads
(from left to right) as follows;
\begin{equation*}
w^{(n)}_{k+d_1+\cdots+d_{n-1}}\leq \cdots \leq
w^{(3)}_{k+d_1+d_2}\leq w^{(2)}_{k+d_1}\leq w^{(1)}_{k}\ \ \ \ \ \
(k\in\Z).
\end{equation*}
Here we assume that $w^{(i)}_k$ is empty if there is no
corresponding entry (see Figure 1).

\begin{figure}\label{semi-infinite ex}
\begin{equation*}
{\bf w} =\ \ \begin{array}{cccc}
 w^{(4)}\!\!  &  w^{(3)}\!\! & w^{(2)}\!\! & w^{(1)}\!\!  \\
  \vdots\!\! &  \vdots\!\! & \vdots\!\! & \vdots\!\!   \\
\!\!-4 & \!\!-4 & \!\!-4 & \!\!-2 \\
\!\! -3 & \!\!-2 & \ 1 &\ 2 \\
\!\! -2 & \!\!-1 & \ 3  &\ 3  \\
\!\! -1 & \ 1 &   &  \\
\ 1 &\ 2 &   &  \\
\ 2 &\ 3 & &   \\
\ 3 &  & & \\
\ 5 &  & &
\end{array} \ , \ \ H_{\Lambda_{\lambda}}=
\begin{array}{cccc}
 w^{(4)}\!\! &  w^{(3)}\!\! & w^{(2)}\!\! & w^{(1)}\!\!  \\
 \vdots\!\! & \vdots\!\! & \vdots\!\! &\vdots\!\!   \\
\!\!-4 & \!\!-4 & \!\!-4 & \!\!-4 \\
\!\! -3 & \!\!-3 &\!\!-3 &\!\!-3 \\
\!\! -2 & \!\!-2 & \!\!-2  &\!\!-2  \\
\!\! -1 & \!\!-1 &   &  \\
\ 1 &\ 1 &   &  \\
\ 2 &\ 2 & &   \\
\ 3 &  & & \\
\ 4 &  & &
\end{array}.
\end{equation*}
\caption{Semi-infinite semistandard tableaux of shape
$\lambda=(4,2,-1,-1)$}
\end{figure}

Let $\Lambda_{\lambda}=\sum_{k=1}^n\Lambda_{\lambda_k}\in P^+$ and
let $\B(\Lambda_{\lambda})$ be the set of all semi-infinite
semistandard tableaux of shape $\lambda$. We may assume that
$\w=w^{(1)}\otimes\cdots\otimes w^{(n)}\in
\B(\Lambda_{\lambda_n})\otimes\cdots\otimes
\B(\Lambda_{\lambda_1})$. By similar arguments as in the case of
usual semistandard tableaux (cf.\cite{KN}), we can check that
$\B(\Lambda_\lambda)$ together with ${\bf 0}$ is stable under
$\te_i$ and $\tf_i$ $(i\in \Z)$. Hence, we have

\begin{prop}\label{BLambdalambda}
$\B(\Lambda_\lambda)$ is a $\gl_{\infty}$-crystal and
$$
\B(\Lambda_\lambda)=
\{\,\tf_{i_1}\cdots\tf_{i_r}H_{\Lambda_{\lambda}}\,|\,r\geq 0, \
i_1,\ldots,i_r\in \Z\,\}\setminus\{{\bf 0}\},
$$
where
$H_{\Lambda_\lambda}=H_{\Lambda_{\lambda_n}}\otimes\cdots\otimes
H_{\Lambda_{\lambda_1}}$. \qed
\end{prop}

Now, let us consider $\B(n\Lambda_0)$ for $n\in\mathbb{N}$. Given
$\w=(w^{(1)},\cdots,w^{(n)})\in \B(n\Lambda_0)$, let $\w_{>0}$ and
$\w_{<0}$ be the subtableaux of $\w$ consisting of positive and
negative entries, respectively. Note that $\w_{>0}\in
SST_{\B_{>0}}(\mu^\pi)$ for some $\mu\in\mathscr{P}$ with $\mu_1\leq
n$, and $\w_{<0}$ is a semi-infinite semistandard tableau of shape
$(-\mu'_n,\ldots,-\mu_1')$.

Suppose that $\w_{<0}=(w^{(1)}_{<0},\ldots,w^{(n)}_{<0})$. For each
$1\leq k \leq n$, we have $w^{(k)}_{<0}\in \B(\Lambda_{-\mu'_{k}})$
and ${\rm wt}(w^{(k)}_{<0})=\Lambda_0-\sum_{i\in I_k}\epsilon_i$ for
a unique $I_k=\{\,-i_{k,1}>\ldots >-i_{k,-\mu'_k}\,\}\subset
\B_{<0}$. We let $\w_{<0}^\vee$ be the tableau of shape $\mu^\pi$,
whose $k$th column (from the right) is filled with
$\{\,-i^\vee_{k,1}<\ldots< -i^\vee_{k,-\mu'_k}\,\}\subset
\B^\vee_{<0}$. It is not difficult to see that $\w_{<0}^\vee\in
SST_{\B^\vee_{< 0}}(\mu^\pi)$. Now, we define
\begin{equation}\label{Psi}
\Psi_n(\w\otimes t_{-n\Lambda_0}) =
\kappa^{-1}(\w^\vee_{<0},\w_{>0})\in\M.
\end{equation}

\begin{ex}{\rm
Let $\w\in\B(4\Lambda_0)$ be as follows.
\begin{equation*}
\begin{array}{cccc}
\vdots\!\! & \vdots\!\! &\vdots\!\! &\vdots\!\!   \\
\!\!-5 & \!\!-5 & \!\!-5 & \!\!-5 \\
\!\!-4 & \!\!-4 & \!\!-3 & \!\!-2 \\
\!\! -3 & \!\!-3 & \!\!-2 & \!\!-1 \\
\!\! -2 & \!\!-1 & \ 1  & \ 2 \\
\ 1 & \ 1 & \ 3  &\  4 \\
\end{array} \ \ .
\end{equation*}
Then
$$
\w^\vee_{<0}=
\begin{array}{cccc}
   &  & -1^\vee\!\!\! & -3^\vee\!\!\! \\
  -1^\vee\!\!\! & -2^\vee\!\!\! & -4^\vee\!\!\! & -4^\vee\!\!\!
\end{array}\ \ , \ \
\w_{>0}=
\begin{array}{cccc}
     &  & 1\ \  & 3\ \  \\
  1\ \  & 1\ \ & 3\ \ & 4\ \
\end{array}\ \ .
$$
Therefore, $$\kappa^{-1}(\w^\vee_{<0},\w_{>0})= \left(
  \begin{array}{cccc}
    1 & 0 & 1 & 0 \\
    0 & 0 & 0 & 1 \\
    1 & 0 & 0 & 0 \\
    1 & 0 & 1 & 0 \\
  \end{array}
\right)\ \ .
$$ }
\end{ex}

\begin{prop}\label{embedding1} For $n\geq 1$, the map
$$\Psi_n : \B(n\Lambda_0)\otimes T_{-n\Lambda_0} \longrightarrow
\M$$
is a $\gl_{\infty}$-crystal embedding.
\end{prop}
\pf (1) Let $\w$ be given. For $i\in\Z_{>0}$ and $x=e,f$, if $\tilde{x}_i\w\neq
\bf{0}$, then
\begin{equation*}
\begin{split}
\Psi_n(\tilde{x}_i(\w\otimes
t_{-n\Lambda_0}))&=\Psi_n((\tilde{x}_i\w)\otimes t_{-n\Lambda_0})
\\
&=\kappa^{-1}(\w^\vee_{<0},\tilde{x}_i\w_{>0}) \\
&=\tilde{x}_i\kappa^{-1}(\w^\vee_{<0},\w_{>0}) \ \ \ \text{by Proposition \ref{bicrystal iso}} \\
&=\tilde{x}_i\Psi_n(\w\otimes t_{-n\Lambda_0}).
\end{split}
\end{equation*}
Similarly, for $i\in\Z_{<0}$ and $x=e,f$, if
$\tilde{x}_i\w\neq \bf{0}$, then
\begin{equation*}
\begin{split}
\Psi_n(\tilde{x}_i(\w\otimes
t_{-n\Lambda_0}))&=\Psi_n((\tilde{x}_i\w)\otimes t_{-n\Lambda_0}) \\
&=\kappa^{-1}(\tilde{x}_i\w^\vee_{<0},\w_{>0})\ \ \ \text{by \cite[Lemma 5.8]{K07}} \\
&=\tilde{x}_i\kappa^{-1}(\w^\vee_{<0},\w_{>0})\ \ \ \text{by Proposition \ref{bicrystal iso}} \\
&=\tilde{x}_i\Psi_n(\w\otimes t_{-n\Lambda_0}).
\end{split}
\end{equation*}
 Finally, comparing the definitions of $\tilde{x}_0$
($x=e,f$) on $\B(n\Lambda_0)$ and $\N$, it is straightforward to see
that $$\Psi_n(\tilde{x}_0(\w\otimes
t_{-n\Lambda_0}))=\kappa^{-1}(\tilde{x}_0(\w^\vee_{<0},\w_{>0})).$$
Since $\kappa$ commutes with $\te_0$ and $\tf_0$ by Theorem \ref{RSK
iso}, we have $\Psi_n(\tilde{x}_0(\w\otimes
t_{-n\Lambda_0}))=\tilde{x}_0\Psi_n(\w\otimes t_{-n\Lambda_0})$. The
other conditions for $\Psi_n$ to be a morphism can be verified
directly. \qed

\begin{rem}\label{imagePsi}{\rm
We have
\begin{equation*}\label{imagePsin}
\begin{split}
&{\rm Im}\Psi_n=\kappa^{-1}\left(\bigsqcup_{\mu \in\mathscr{P},\
\mu_1\leq n}SST_{\B^\vee_{<0}}(\mu^\pi)\times
SST_{\B_{>0}}(\mu^\pi)\right), \\
& {\rm Im}\Psi_n \subset {\rm Im}\Psi_{n+1} \ \ (n\geq 1),\\
& \M=\bigcup_{n\geq 1}{\rm Im}\Psi_n.
\end{split}
\end{equation*}
Note that there exists a strict morphism $\Phi_n :
\M\otimes T_{ n\Lambda_0} \longrightarrow \B(n\Lambda_0)$ sending
$\mathbb{O}\otimes t_{n\Lambda_0}$ to $H_{n\Lambda_0}$ such that
$\Phi_n(A)\neq {\bf 0}$ if and only if $A\in {\rm Im}\Psi_n$. }
\end{rem}

\subsection{Crystal graphs of generalized Verma modules}

Given  $\mu,\nu\in \cP$, we put
\begin{equation}
\begin{split}
\M_{\mu,\nu}&=\M \times SST_{\B^\vee_{<0}}(\mu^\pi)\times
SST_{\B_{>0}}(\nu), \\
\mathbb{O}_{\mu,\nu}&=(\mathbb{O},H_{\mu^\pi},H_{\nu}).
\end{split}
\end{equation}
For $(A,S_{<0},S_{>0})\in\M_{\mu,\nu}$, we define
\begin{itemize}
\item[(1)] if $\tilde{x}_i(A \otimes S_{<0} )=A'\otimes S_{<0}'$ for
$i\in\Z_{<0}$, then
$$\tilde{x}_i(A,S_{<0},S_{>0})=(A',S_{<0}',S_{>0}),$$

\item[(2)] if $\tilde{x}_i(A\otimes S_{>0})=A''\otimes S_{>0}''$ for
$i\in\Z_{>0}$, then
$$\tilde{x}_i(A,S_{<0},S_{>0})=(A'',S_{<0},S_{>0}''),$$

\item[(3)] $\tilde{x}_0(A,S_{<0},S_{>0})=(\tilde{x}_0A,S_{<0},S_{>0})$,
\end{itemize}
where $x=e,f$, and $\widetilde{x}_i(A,S_{<0},S_{>0})={\bf 0}$ if any
of its components is ${\bf 0}$.

\begin{prop}
For $\mu,\nu\in \cP$, $\M_{\mu,\nu}$ is a $\gl_{\infty}$-crystal
 and
$$
\M_{\mu,\nu}=
\{\,\tf_{i_1}\cdots\tf_{i_r}\mathbb{O}_{\mu,\nu}\,|\,r\geq
0, \ i_1,\ldots,i_r\in \Z\,\}\setminus\{{\bf 0}\}.
$$
\end{prop}
\pf Let $(A,S_{<0},S_{>0})\in\M_{\mu,\nu}$ be given with
$\te_i(A,S_{<0},S_{>0})={\bf 0}$ for all $i\in\Z$. First, we have
$\te_iA={\bf 0}$ for all $i\in\Z^\times$, which implies that $A$ is
a diagonal matrix with entries $a_{-1^\vee, 1}\geq a_{-2^\vee,
2}\geq \ldots$. Since $\te_0A={\bf 0}$, we have $a_{-1^\vee, 1}=0$
and hence $A=\mathbb{O}$. This implies that $S_{<0}=H_{\mu^\pi}$ and
$S_{>0}=H_{\nu}$ since $\te_i S_{<0}={\bf 0}$ for $i\in\Z_{<0}$ and
$\te_iS_{>0}={\bf 0}$ for $i\in\Z_{>0}$, respectively. \qed

For $n\geq \mu_1+\nu_1$, we put $$\Lambda_{\mu,\nu\, ;\,
n}=\Lambda_{\lambda},$$ where $\lambda=(\lambda_1,\ldots,\lambda_n)$
is a generalized partition of length $n$ such that
\begin{equation*}\label{lambda}
\begin{split}
\lambda^+&=(\max(\lambda_1,0),\ldots,\max(\lambda_n,0))=\nu', \\
\lambda^-&=(\max(-\lambda_n,0),\ldots,\max(-\lambda_1,0))=\mu'.
\end{split}
\end{equation*}

\begin{prop}\label{embedding}
Let $\mu,\nu\in \cP$ be given. Then for $n\geq \mu_1+\nu_1$, there
exists an embedding of $\gl_{\infty}$-crystals
$$\Psi_{\mu,\nu\,;\,n} : \B(\Lambda_{\mu,\nu\,;\,n})\otimes T_{-n\Lambda_0} \longrightarrow \M_{\mu,\nu},$$
sending $H_{\Lambda_{\mu,\nu;n}}\otimes t_{-n\Lambda_0}$ to
$\mathbb{O}_{\mu,\nu}$.
\end{prop}
\pf Given $\w=(w^{(1)},\cdots,w^{(n)})\in
\B(\Lambda_{\mu,\nu\,;\,n})$, consider the subtableau of $\w$
consisting of positive entries, say $\w_{>0}$. Let $\w^+_{>0}$ be
the subtableau of $\w_{>0}$ corresponding to the positions where
$H_{\Lambda_{\mu,\nu\,;\,n}}$ has positive entries, and let
$\w^-_{>0}$ be the compliment of $\w^+_{>0}$ in $\w_{>0}$. Note that
$\w^+_{>0}\in SST_{\B_{>0}}(\nu)$ and $\w^-_{>0}\in
SST_{\B_{<0}}((\eta/\mu)^\pi)$ for some $\eta \supset \mu$ with
$\eta'_1\leq n$.

Let $\w_{<0}$ be the subtableau of $\w$ consisting of negative
entries, which is also a semi-infinite semistandard tableau of shape
$(-\eta'_n,\ldots,-\eta_1')$. By the same method as in
(\ref{BLambdalambda}), we obtain $\w_{<0}^\vee\in
SST_{\B^\vee_{<0}}(\eta^\pi)$ from $\w_{<0}$.

Given $A=A(\bi,\bj)\in \M$ and $S\in SST_{\B^\vee_{<0}}(\mu^\pi)$,
we define ${\bf P}(S \leftarrow A)$ to be the tableau obtained by
inserting $\bi_{\rm rev}=i_r\cdots i_1$ to $S$, that is,
$${\bf P}(S \leftarrow A)=(\cdots((S\leftarrow{i_r}\,)\leftarrow{i_{r-1}}\,)
\cdots)\leftarrow{i_1}.$$ Suppose that ${\rm sh}{\bf P}(S\leftarrow
A)=\tau^\pi$ for some $\tau\supset \mu$. For $1\le k \le r$, let us
fill the box in $(\tau/\mu)^\pi$ with $c$ if it is created when
$i_k$ is inserted into $(({S}\leftarrow{i_{r}}\,)
\cdots)\leftarrow{i_{k+1}}$ and $j_k=c$. This defines the recording
tableau ${\bf Q}(S\leftarrow A)$ of shape $(\tau/\mu)^\pi$.

Now, we let $A\in \M$ and $S_{<0}\in SST_{\B^\vee_{<0}}(\mu^\pi)$ be
the unique pair such that ${\bf P}(S_{<0}\leftarrow A)=\w_{<0}^\vee$
and ${\bf Q}(S_{<0}\leftarrow A)=\w_{>0}^-$, and let
$S_{>0}=\w_{>0}^+$. This defines a map $\Psi_{\mu,\nu\,;\,n} :
\B(\Lambda_{\mu,\nu\,;\,n})\otimes T_{-n\Lambda_0} \rightarrow
\M_{\mu,\nu}$ by $$\Psi_{\mu,\nu\,;\,n}(\w\otimes
t_{-n\Lambda_0})=(A,S_{<0},S_{>0}).$$ Modifying the arguments in
Theorem \ref{RSK iso} and Proposition \ref{embedding1}, we can check
that $\Psi_{\mu,\nu\,;\,n}$ is a $\gl_{\infty}$-crystal embedding. \qed

\begin{ex}{\rm
Let $\w\in\B(\Lambda_{(4,2,-1,-1)})$ be as in Figure 1. Note that
$\mu=(2)$,  $\nu=(2,2,1,1)$ and $n=4$. We have
\begin{equation*}
\begin{split}
\w^\vee_{<0}=
\begin{array}{ccc}
    & -1^\vee\!\!\! & -1^\vee\!\!\!  \\
    & -2^\vee\!\!\! & -3^\vee\!\!\! \\
   -3^\vee\!\!\! & -3^\vee\!\!\! & -4^\vee\!\!\!
\end{array}\ , \ \
\w^-_{>0}=
\begin{array}{ccc}
   & 1\ & 2\ \\
   & 3\  & 3\  \\
  1\ & \Box\  & \Box\ \\
\end{array}\ , \ \ \w^+_{>0}=
\begin{array}{cc}
 1\ & 2\   \\
 2\ & 3\    \\
 3\ &   \\
 5\ &
\end{array} \ \ .
\end{split}
\end{equation*}
Then we can check that ${\bf P}(S \leftarrow A)=\w^\vee_{<0}$ and
${\bf Q}(S \leftarrow A)=\w^-_{>0}$, where
\begin{equation*}
\begin{split}
S=
\begin{array}{ccc}
   &\ \ \Box\  &\ \ \Box \\
   &\ \ \Box\ & \ \ \Box \\
 \Box  & -3^\vee\!\!\!  & -4^\vee\!\!\! \\
\end{array} \ \ , \ \ A= \left(
  \begin{array}{ccc}
    1 & 1 & 0 \\
    0 & 0 & 1 \\
    1 & 0 & 1 \\
  \end{array}
\right)\ \ .
\end{split}
\end{equation*}
Hence
\begin{equation*}
\Psi_{(2),(2,2,1,1)\,;\,4}(\w\otimes t_{-4\Lambda_0})=\left(
\left(\begin{array}{ccc}
    1 & 1 & 0 \\
    0 & 0 & 1 \\
    1 & 0 & 1 \\
  \end{array}\right)\ , \ -3^\vee\!\!
-4^\vee\ ,\
\begin{array}{cc}
 1\ & 2   \\
 2\ & 3    \\
 3\ &   \\
 5\ &
\end{array}\right) \in \M_{(2),(2,2,1,1)}.
\end{equation*}
 }
\end{ex}

\begin{rem}{\rm
(1) In case of $\nu=\emptyset$, the map sending
$(\w^\vee_{<0},\w^-_{>0})$ to $(A,S_{<0})$ is a skew version of the
RSK correspondence introduced by Sagan and Stanley \cite{SS}.

(2) Note that $\M_{\mu,\nu}=\bigcup_{n\geq \mu_1+\nu_1}{\rm
Im}\Psi_{\mu,\nu\,;\,n}$. Hence, we may regard $\M_{\mu,\nu}$ as the
limit of $\B(\Lambda_{\mu,\nu\,;\,n})$ as $n\rightarrow \infty$ and
the crystal graph of the generalized Verma module induced from an
irreducible $\frak{l}$-module with $\frak{l}$-dominant highest
weight $-\sum_{i<0}\mu_{-i}\epsilon_i+\sum_{j>0}\nu_j\epsilon_j$. We
also have a strict morphism $\Phi_{\mu,\nu\,;\,n} :
\M_{\mu,\nu}\otimes T_{ n\Lambda_0} \rightarrow
\B(\Lambda_{\mu,\nu\,;\,n})$ such that
$\Phi_{\mu,\nu\,;\,n}((A,S_{<0},S_{>0})\otimes t_{ n\Lambda_0})=\w
\neq \bf{0}$ if and only if $\Psi_{\mu,\nu\,;\,n}(\w\otimes
t_{-n\Lambda_0})=(A,S_{<0},S_{>0})$.

(3) Bitableaux realizations of irreducible highest weight
representations of Lie (super) algebras including
$\B(\Lambda_{\lambda})$ and their combinatorics  can be found in \cite{K08}. }
\end{rem}

\section{Demazure crystals and a flagged RSK correspondence}

\subsection{Demazure crystals}\mbox{}
For $i\in\mathbb{Z}$, let $s_i\in GL(\frak{h}^*)$ be the simple
reflection with respect to $\alpha_i$ defined by
$$s_i(\lambda)=\lambda-\langle\lambda,h_i\rangle\alpha_i, \ \ (\lambda\in\frak{h}^*). $$
Then $s_i$ acts as the transposition on
$\{\,\epsilon_i\,|\,i\in\mathbb{Z}^{\times}\,\}$ (hence on
$\mathbb{Z}^{\times}$) given by
\begin{equation*}
s_i=
\begin{cases}
(i\ i+1), & \text{if $i>0$}, \\
(i-1\ i), & \text{if $i<0$}, \\
(-1\ 1), & \text{if $i=0$}.
\end{cases}
\end{equation*}
Let $W$ be the Weyl group of $\gl_{\infty}$, which is generated by
$\{\,s_i\,|\,i\in\mathbb{Z}\,\}$, and let $\ell(w)$ denote the length of
$w\in W$. For $\Lambda\in P^+$, let $W_{\Lambda}$ be the stabilizer
of $\Lambda$, and let $W^{\Lambda}=\{\,w\,|\,\ell(ws_i)>\ell(w)
\text{ for $s_i\in W_{\Lambda}$}\,\}$.   Let $<$ denote the Bruhat
order on $W$. It induces the Bruhat order on $W^{\Lambda}$, which is
also denoted by $<$.

Let $w\in W^{\Lambda}$ be given and $w=s_{i_1}\cdots s_{i_r}$ its
reduced expression. We define the {\it Demazure crystal of
$\B(\Lambda)$ associated with $w$} \cite{Kas93} by
\begin{equation}
\B_w(\Lambda)=\{\,\tf_{i_1}^{m_1}\cdots\tf_{i_r}^{m_r}H_{\Lambda}\,|\,m_1,\ldots,m_r\geq
0\,\}\setminus \{{\bf 0}\}.
\end{equation}
For any $i\in \Z$, an $i$-string $S$ in $\B(\Lambda)$, a connected
component with only $i$-arrows, satisfies one of the following three
conditions;
\begin{equation}\label{Demazure}
\begin{split}
& S \subset \B_w(\Lambda), \\
& S\cap \B_w(\Lambda)=\emptyset, \\
& S\cap \B_w(\Lambda) \text{ is a highest weight vector of $S$}.
\end{split}
\end{equation}
Now, given $\mu,\nu\in\cP$, we put
\begin{equation}\label{Mw}
\M_{\mu,\nu,w}=\{\,\tf_{i_1}^{m_1}\cdots\tf_{i_r}^{m_r}\mathbb{O}_{\mu,\nu}\,|\,m_1,\ldots,m_r\geq
0\,\}\setminus \{{\bf 0}\}.
\end{equation}
For $w,w'\in W^{\Lambda}$, we have $\M_{\mu,\nu,w}\subset
\M_{\mu,\nu,w'}$ if and only if $w\leq w'$. Since each element in
$\M_{\mu,\nu,w}$ is contained in
$\Psi_{\mu,\nu\,;\,n}(\B_w(\Lambda_{\mu,\nu\,;\,n})\otimes
T_{-n\Lambda_0})$ for a sufficiently large $n$, it is well-defined.


\subsection{Grassmannian permutations}

Let $\lambda$ be a partition. The {\it residue} of a box in
$\lambda$ is given by $j-i$ if it is located in the $i$th row and
the $j$th column. A {\it standard tableau of shape $\lambda$} is a
tableau obtained by filling $\lambda$ with
$\{\,1,\ldots,|\lambda|\,\}$ in such a way that the entries in each
column (resp. row) are increasing from top to bottom (resp. left to
right).

Consider $W_{\Lambda_0}=\langle\,
s_i\,|\,i\in\mathbb{Z}^\times\,\rangle$ and let $w\in W^{\Lambda_0}$
be given. Let $D(w)=\{\,(i,j)\in\Z^{\times}\times
\Z^{\times}\,|\,i<w^{-1}(j),\ j<w(i)\,\}$ be the {\it diagram of
$w$} and let $\lambda(w)=(\lambda(w)_i)_{i\geq 1}$ be the {\it shape
of $w$}, where $\lambda(w)_i=|\{\,j\,|\,(-i,j)\in D(w)\,\}|$. Since
$w(i)<w(i+1)$ for $i\leq -2$ or $i\geq 1$, and $w(-1)> w(1)$,
$\lambda(w)$ is a partition. Conversely, a partition $\lambda$
determines a unique permutation $w\in W^{\Lambda_0}$ such that
$\lambda(w)=\lambda$. Hence, we have a bijection from
$W^{\Lambda_0}$ to $\cP$ sending $w$ to $\lambda(w)$, where
$|\lambda(w)|=\ell(w)$. If $T$ is a standard tableau of shape
$\lambda(w)$ and $a_i$ is the residue of the box corresponding to
$i$ in $T$ ($1\leq i\leq \ell(w)$), then
$w=s_{a_{\ell(w)}}s_{a_{\ell(w)-1}}\cdots s_{a_1}$ gives a reduced
expression of $w$. For $w,w'\in W^{\Lambda_0}$, we have $w\leq w'$
if and only if $\lambda(w)\subseteq \lambda(w')$ (cf.\cite{Mac}).

\begin{ex}{\rm
Let $w\in W^{\Lambda_0}$ be given by
$$
w= \left[
  \begin{array}{cccccccccccccc}
    \cdots &    -5 & -4 & -3 & -2 & -1 & \,\,\,\,1 & \,\,\,\,2 & \,\,\,\,3 & \,\,\,4 & \,\,\,5 & \,\,\,6 &  \cdots \\
    \cdots &    -4 & -2 & \,\,\,\,2 & \,\,\,\,5 & \,\,\,\,6 & -5 & -3 & -1 & \,\,\,1 & \,\,\,3 & \,\,\,4 &  \cdots \\
  \end{array}
\right]\,
$$
where $w(i)=i$ for $i\geq 7$ or $i\leq -6$. Then
$\lambda(w)=(6,6,4,2,1)$, where the residue on each box is given by
$$
\begin{array}{ccccccc}
  0 & 1 & 2 & 3 & 4 & 5 &   \\
 \!\!\!\! -1 & 0 & 1 & 2 & 3 & 4 &   \\
 \!\!\!\! -2 & \!\!\!\! -1 & 0 & 1 &  &  & \\
 \!\!\!\! -3 & \!\!\!\! -2 &  &  &  &  & \\
 \!\!\!\! -4  &  &  &  &  &  & \\
\end{array}\ \ .
$$

 }
\end{ex}

\subsection{Flagged skew Schur functions}
Let $X=\{\,x_1,x_2,\ldots\,\}$ be the set of variables and
$X_k=\{\,x_1,\ldots,x_k\,\}$ for $k\geq 1$. Let
$\phi=(\phi_1,\ldots,\phi_d)$ be a sequence of weakly increasing
positive integers of length $d$, which is called a {\it flag of
length $d$}. Given a skew Young diagram $\lambda/\mu$ with
$\ell(\lambda),\ell(\mu)\leq d$, we define the {\it flagged Schur
function} $s_{\lambda/\mu}(X_{\phi})$ by
\begin{equation}\label{multiSchur}
s_{\lambda/\mu}(X_{\phi})={\rm det}\left(
h_{\lambda_i-\mu_j-i+j}(X_{\phi_i}) \right)_{1\leq i,j\leq d},
\end{equation}
where $h_k(X_{\phi_i})$ is the $k$th complete symmetric function in
$X_{\phi_i}$ (cf.\cite{Mac}). An equivalent definition is that
\begin{equation}
s_{\lambda/\mu}(X_{\phi})=\sum_{T}x^T,
\end{equation}
where the sum ranges over all semistandard tableaux of shape
$\lambda/\mu$ such that the entries in the $i$th row are no more
than $\phi_i$ for $1\leq i\leq d$. Here, $x^T=\prod_{i}x_i^{m_i}$,
where $m_i$ is the number of occurrences of $i$ in $T$.

Let $SST_{\B_{<0}^\vee}(\lambda/\mu)_{\phi}$ (resp.
$SST_{\B_{>0}}(\lambda/\mu)_{\phi}$)   be the set of semistandard
tableaux of shape
$\lambda/\mu$ such that the entries in the $i$th row are no more than
$-\phi_i^\vee$ (resp. $\phi_i$) for $1\leq i\leq d$. Let
$Y=\{\,y_1,y_2,\ldots\,\}$ be another set of variables and
$Y_k=\{y_1,\ldots,y_k\}$ for $k\geq 1$. Then we have
\begin{equation*}
{\rm
ch}SST_{\B_{<0}^\vee}(\lambda/\mu)_{\phi}=s_{\lambda/\mu}(X_{\phi}),
\ \ \ {\rm
ch}SST_{\B_{>0}}(\lambda/\mu)_{\phi}=s_{\lambda/\mu}(Y_{\phi}),
\end{equation*}
where we put $$x_i=e^{{\rm wt}(-i^\vee)}=e^{-\epsilon_{-i}},\ \
y_j=e^{{\rm wt}(j)}=e^{\epsilon_j}$$ for $-i^\vee\in\B^\vee_{<0}$
and $j\in\B_{>0}$. When $S\in SST_{\B^\vee_{<0}}(\nu^\pi)_{\phi}$
(resp. $S\in SST_{\B_{>0}}(\nu^\pi)_{\phi}$) is given for
$\nu\in\cP$ with $\ell(\nu)\leq d$, we   view
$\nu^\pi=(n^d)/(n-\nu_n,\ldots,n-\nu_1)$ for some $n$, and
understand that the entries of $S$ in each $i$th row from the bottom
are no more than $-\phi^\vee_{d-i+1}$ (resp. $\phi_{d-i+1}$) for
$1\leq i\leq d$.

Let $\alpha=(\alpha_1,\ldots,\alpha_d)$ be the sequence of weakly decreasing positive integers of length $d$. Let
\begin{equation}
\widehat{s}_{\lambda/\mu}(X_{\alpha})={\rm det}\left(
h_{\lambda_i-\mu_j-i+j}(X_{\alpha_j}) \right)_{1\leq i,j\leq d}.
\end{equation}
By (\ref{multiSchur}), it is easy to check that
$\widehat{s}_{\lambda/\mu}(X_{\alpha})=s_{\widehat{\mu}/\widehat{\lambda}}(X_{\phi})$,
where $\widehat{\lambda}=(n-\lambda_{d-i+1})_{1\leq i\leq d}$,
$\widehat{\mu}=(n-\mu_{d-i+1})_{1\leq i\leq d}$ for some $n\geq
\lambda_1,\mu_1$, and $\phi$ is the reverse sequence of $\alpha$. In
particular, we have for $\nu\in\cP$ with $\ell(\nu)\leq d$,
\begin{equation*}
s_{\nu^\pi}(X_{\phi})=\widehat{s}_{\nu}(X_{\alpha}).
\end{equation*}

\subsection{A flagged RSK correspondence}
In the sequel, we assume that $S$ is a finite subset of
$\mathbb{N}^2=\mathbb{N}\times \mathbb{N}$. We define $\theta(S)$ to
be the border strip of
the smallest partition $\lambda$ such that $S\subset\lambda$. 
Recall that a {\it border strip} of a partition $\lambda$ is a skew
diagram $\lambda/\mu$, where
$\lambda=(\alpha_1,\ldots,\alpha_d|\beta_1,\ldots,\beta_d)$ and
$\mu=(\alpha_2,\ldots,a_d|\beta_2,\ldots,\beta_d)$ following
Frobenius notation. We put $c(S)=(\alpha_1,\beta_1)$.

\begin{ex}\label{lambda of S}{\rm Let
$S=\{\,(1,1),(1,4),(2,2),(3,1),(3,3),(4,3)\,\}$. Then
$$S\ =\
\begin{tabular}{c|ccccc}
     & 1 & 2 & 3 & 4 & $\cdots$ \\   \hline
  $1$ & $\bullet$ &  &  & $\bullet$ \\
  $2$ &  & $\bullet$ &   &  \\
  $3$ & $\bullet$ &  & $\bullet$ &  \\
  $4$ &  &   & $\bullet$ &  \\
  \vdots &  &  &  &
\end{tabular} \ \ , \ \
\theta(S)=(4,3,3,3)/(2,2,2)=
\begin{array}{cccc}
  \square & \square & \blacksquare & \blacksquare \\
  \square & \square & \blacksquare &   \\
  \square & \square & \blacksquare &   \\
  \blacksquare & \blacksquare & \blacksquare  & \ast
\end{array}\ \ ,
$$
where the skew diagram consisting of the black boxes is the border
strip $\theta(S)$ and $\ast$ indicates the point $c(S)=(4,4)$.}
\end{ex}

We define inductively a finite sequence of points
$c_1=(\alpha_1,\beta_1),\ldots,c_d=(\alpha_d,\beta_d)$ as follows;

\begin{itemize}
\item[(1)] let $S^{(1)}=S$ and put $c_1=c(S^{(1)})$,

\item[(2)] for $1\leq k\leq d-1$, let $S^{(k+1)}=S^{(k)}\setminus \theta(S^{(k)})$,
and put $c_{k+1}=c(S^{(k+1)})$,
\end{itemize}
where $d$ is the smallest one such that $S^{(d+1)}=\emptyset$. Note
that $c_{k+1}$ is located to the northwest of $c_k$, that is,
$\alpha_k> \alpha_{k+1}$ and $\beta_k > \beta_{k+1}$. Now, we define
\begin{equation}
\lambda(S)=(\alpha_1,\ldots,\alpha_d | \beta_1,\ldots,\beta_d)\in
\cP
\end{equation}
following Frobenius notation, where $(\alpha_1,\ldots,\alpha_d)$
(resp. (resp. $(\beta_1,\ldots,\beta_d)$) corresponds to the lower
(resp. upper) flag of $\lambda(S)$ with respect to its diagonal. We
define $w(S)$ to be the permutation in $W^{\Lambda_0}$ corresponding
to $\lambda(S)$, that is,
\begin{equation}
\lambda(w(S))=\lambda(S).
\end{equation}

\begin{ex}{\rm Let $S$ be as in
Example \ref{lambda of S}. Then {\allowdisplaybreaks
\begin{equation*}
\begin{split}
& S^{(1)}\ = \
\begin{tabular}{c|cccc}
     & 1 & 2 & 3 & 4  \\   \hline
  $1$ & $\bullet$ &  &  & $\bullet$ \\
  $2$ &  & $\bullet$ &   &  \\
  $3$ & $\bullet$ &  & $\bullet$ &  \\
  $4$ &  &  & $\bullet$ &  \\
\end{tabular}\ \ \ \ , \ \ \
\theta(S^{(1)})=(4,3,3,3)/(2,2,2)=
\begin{array}{cccc}
  \square & \square & \blacksquare & \blacksquare \\
  \square & \square & \blacksquare &   \\
  \square & \square & \blacksquare &   \\
  \blacksquare & \blacksquare & \blacksquare  & \ast
\end{array}\ \ ,\\
\end{split}
\end{equation*}\vskip 2mm
\begin{equation*}
\begin{split}
& S^{(2)}\ = \
\begin{tabular}{c|cccc}
     & 1 & 2 & 3 & 4  \\   \hline
  $1$ & $\bullet$ &  &  &   \\
  $2$ &  & $\bullet$ &   &  \\
  $3$ & $\bullet$ &  &   &  \\
  $4$ &  &   &  &  \\
\end{tabular}\ \ \ \ \ , \ \ \ \ \ \ \ \ \ \ \
\theta(S^{(2)})=(2,2,1)/(1)=
\begin{array}{cccc}
  \square & \blacksquare & \cdot & \cdot \\
  \blacksquare & \blacksquare & \cdot &   \\
  \blacksquare & \ast & \cdot &   \\
  \cdot & \cdot &   &
\end{array}\ \ \ \  \ ,\\
\end{split}
\end{equation*}\vskip 2mm
\begin{equation*}
\begin{split}
& S^{(3)}\ =\
\begin{tabular}{c|cccc}
     & 1 & 2 & 3 & 4  \\   \hline
  $1$ & $\bullet$ &  &  &   \\
  $2$ &  &    &   &  \\
  $3$ &   &  &   &  \\
  $4$ &  &   &  &  \\
\end{tabular}\ \ \ \ \  \ , \ \ \ \ \ \ \ \ \ \ \ \ \ \ \ \ \ \,
\
\theta(S^{(3)})\ =\ (1)=
\begin{array}{cccc}
  \blacksquare & \cdot &   &   \\
  \cdot & \cdot &   &   \\
  \cdot &   &   &   \\
    &   &   &
\end{array}\ \ \ \ \ \ \ \ \ ,
\end{split}
\end{equation*}}
where we have $c_1=(4,4)$, $c_2=(3,2)$,
$c_3=(1,1)$. Hence
$$\lambda(S)=(4,3,1|4,2,1)=(4,4,3,2).$$ }
\end{ex}

For $w\in W^{\Lambda_0}$ with a reduced expression $w=s_{i_1}\cdots
s_{i_r}$, we put
\begin{equation}
\M_w=\{\,\tf_{i_1}^{m_1}\cdots\tf_{i_r}^{m_r}\mathbb{O}\,|\,m_1,\ldots,m_r\geq
0\,\}\setminus \{{\bf 0}\},
\end{equation}
that is, $\M_w=\M_{\emptyset,\emptyset,w}$ (see (\ref{Mw})). For
$A\in \M$, let ${\rm supp}(A)=\{\,(i,j)\,|\,a_{-i^\vee, j}\neq
0\,\}\subset \mathbb{N}^2$ be the support of $A$.

\begin{thm}\label{Demazure-0} For $w\in W^{\Lambda_0}$, we have
\begin{equation*}
\M_w=\{\,A\,|\ \lambda({\rm supp}(A))\subset \lambda(w)\
\}=\{\,A\,|\ w({\rm supp}(A))\leq w\ \}.
\end{equation*}
\end{thm}
\pf It suffices to prove the first identity. We use induction on
$\ell(w)=|\lambda(w)|$.

If $|\lambda(w)|=1$, then $w=s_0$ and it is clear. We assume that
$|\lambda(w)|\geq 2$. Choose $w'\in W^{\Lambda_0}$ such that
$\ell(w')=\ell(w)-1$, equivalently $\lambda(w')\subset \lambda(w)$
with $|\lambda(w)/\lambda(w')|=1$. Let $r$ be the residue of
$\lambda(w)/\lambda(w')$.

Choose a standard tableau $T$ of shape $\lambda(w)$ such that the
largest entry occurs at $\lambda(w)/\lambda(w')$. Let $a_i$ be the
residue of the box corresponding to $i$ in $T$ ($1\leq i\leq
\ell(w)$). Then we have reduced expressions
$w=s_{a_{\ell(w)}}s_{a_{\ell(w)-1}}\cdots s_{a_1}$  and
$w'=s_{a_{\ell(w)-1}}\cdots s_{a_1}$, where $a_{\ell(w)}=r$. Note
that
\begin{equation}\label{inductive relation}
\M_w=\bigcup_{k\geq 0}\tf_r^k\M_{w'} \setminus \{{\bf 0}\}.
\end{equation}
Let $\mathcal{N}_w=\{\,A\,|\ \lambda({\rm supp}(A))\subset
\lambda(w)\ \}$. By induction hypothesis, it suffices to show that
$\mathcal{N}_w=\bigcup_{k\geq 0}\tf_r^k\mathcal{N}_{w'} \setminus
\{{\bf 0}\}$.

Let $A\in \mathcal{N}_{w}$ be given. We first claim that $A\in\tf_r^k
\mathcal{N}_{w'}$ for some $k\geq 0$. We will keep the previous
notations $\lambda(S)$, $\theta(S^{(k)})$, and
$c_k=(\alpha_k,\beta_k)$ ($1\leq k\leq d$) with $S={\rm supp}(A)$.
If $\lambda(S)$ does not contain the box $c$ corresponding to
$\lambda(w)/\lambda(w')$, then $\lambda(S) \subset \lambda(w')$ and
$A\in \mathcal{N}_{w'}$. So we may assume that $c\in \lambda(S)$.

\textsc{Case 1}. Suppose that $r=0$. By definition of $\lambda(S)$,
we have $\theta(S^{(d)})=(1,1)$. If we choose $k\geq 1$ such that
the entry of $\te^k_0A$ at $(1,1)$ is $0$, then ${\rm
supp}(\te^k_0A)={\rm supp}(A)\setminus\{(1,1)\}$. Hence
$\lambda({\rm supp}(\te^k_0A))\subset \lambda(w')$ and
$\te^k_0A\in\mathcal{N}_{w'}$.

\textsc{Case 2}. Suppose that $r\neq 0$. We may assume that $r>0$
since the argument for $r<0$ is almost the same. In this case, we
have $c_s=(\alpha_s,\beta_s)$ with $\beta_s=r+1$ for some $1\leq
s\leq d$. There exists $1\leq m\leq \alpha_s$ such that
$a_{-m^\vee,r+1}\neq 0$ with $(m,r+1)\in S^{(s)}$ and
$(i,r+1)\not\in S^{(s)}$ for $m<i\leq \alpha_s$. Since $c$ is a
removable corner of $\lambda(w)$ and hence of $\lambda(S)$, we have
$a_{-i^\vee, r}=0$ for $1\leq i < m$ (see Figure 2). Otherwise, we
have $\beta_{s-1}=r$, which implies that $c$ is not removable. Let
$A=A(\bi,\bj)$ and consider the subword of $\bi$ consisting of $r$
and $r+1$. Then the $r+1$'s corresponding to the non-zero entries
$a_{-i^\vee, r+1}$ for $1\leq i \leq m$ can be replaced by $r$
applying $\te_r^k$ for some $k\geq 1$. Equivalently, applying
$\te_r^k$, the entries $a_{-i^\vee, r}$ (resp. $a_{-i^\vee, r+1}$)
are replaced by $a_{-i^\vee, r}+a_{-i^\vee, r+1}$ (resp. $0$) for
$1\leq i\leq m$. On the other hand, we can check that
$(\alpha_t,\beta_t)$ are invariant under $\te_r^k$ for $1\leq t\leq
d$ with $t\neq s$. Hence $\lambda({\rm supp}(\te_r^k A))\subset
\lambda(w')$ and $\te_r^k A\in\mathcal{N}_{w'}$.

Conversely, let $A\in \mathcal{N}_{w'}$ be given. Using similar
arguments, it is not difficult to check that either $\lambda({\rm
supp}(\tf^k_r A))=\lambda({\rm supp}(A))$ or $\lambda({\rm
supp}(\tf^k_r A))/\lambda({\rm supp}(A))=c$ whenever $\tf^k_r A\neq
{\bf 0}$. Hence   $\lambda({\rm supp}(\tf^k_r A))\subset
\lambda(w)$. This completes our induction. \qed

\begin{figure}
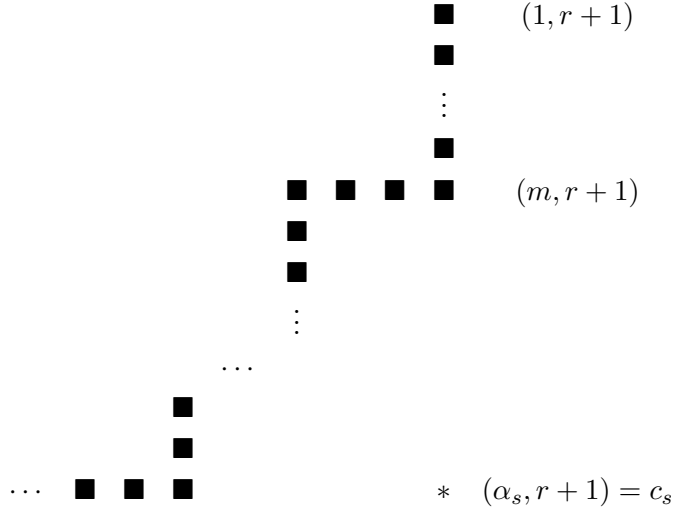

$$
\begin{array}{cccccccccccc}
    &   &   &   &   &   &   &   &   &   & \blacksquare  & (1,r+1) \\
    &   &   &   &   &   &   &   &   &   & \blacksquare  &   \\
    &   &   &   &   &   &   &   &   &   & \vdots  &   \\
    &   &   &   &   &   &   &   &   &   & \blacksquare  &   \\
    &   &     &  &    &   &   & \blacksquare  & \blacksquare  & \blacksquare  & \blacksquare  & (m,r+1)  \\
    &   &   &   &   &   &   & \blacksquare  &   &   &   &   \\
    &   &   &   &   &   &   & \blacksquare  &   &   &   &   \\
    &   &   &   &   &   &   & \vdots  &   &   &   &   \\
    &   &   &   &   &   & \cdots  &   &   &   &   &   \\
    &   &   &   &   & \blacksquare  &   &   &   &   &   &   \\
    &   &   &   &   & \blacksquare  &   &   &   &   &   &   \\
    &   &  \cdots & \blacksquare  &  \blacksquare &  \blacksquare &   &   &   &   &  \ast &  (\alpha_s,r+1)=c_s  \\
\end{array}
$$
\caption{The border strip $\theta(S^{(s)})$}
\end{figure}

\begin{cor}
For $w\in W^{\Lambda_0}$, we have
\begin{equation*}
{\rm ch}\M_w=\sum_{\substack{S\subset\mathbb{N}^2 \\
w(S)\leq w}}\prod_{(i,j)\in S}\frac{x_iy_j}{(1-x_iy_j)}.
\end{equation*}
\end{cor}

\begin{rem}\label{rootexpression}{\rm
Identifying $(i,j)\in \mathbb{N}^2$  with
$-\epsilon_{-i}+\epsilon_j\in\Delta(\mathfrak{u}_-)$, one may write
\begin{equation*}
{\rm ch}\M_w=\sum_{\substack{S\subset\Delta(\mathfrak{u}_-) \\
w(S)\leq w}}\prod_{\alpha\in S}\frac{e^\alpha}{(1-e^\alpha)}.
\end{equation*}}
\end{rem}

\vskip 5mm Now, let us put
\begin{equation}
\N_w=\kappa(\M_w)=\{\,\tf_{i_1}^{m_1}\cdots\tf_{i_r}^{m_r}(\emptyset,\emptyset)\,|\,m_1,\ldots,m_r\geq
0\,\}\setminus \{{\bf 0}\}.
\end{equation}
for $w\in W^{\Lambda_0}$ with a reduced expression $w=s_{i_1}\cdots
s_{i_r}$. We define $\alpha(w)=(\alpha_1,\ldots,\alpha_d)$ and
$\beta(w)=(\beta_1,\ldots,\beta_d)$ to be strict partitions of length
$d$ such that
\begin{equation}\label{lambda w}
\lambda(w)=(\alpha(w)|\beta(w)).
\end{equation}
We put $d(w)=d$, the diagonal length of $\lambda(w)$ and
\begin{equation}
\begin{split}
\phi(w)&=(\phi_1,\ldots,\phi_d)=(\alpha_d,\ldots,\alpha_1), \\
\psi(w)&=(\psi_1,\ldots,\psi_d)=(\beta_d,\ldots,\beta_1),
\end{split}
\end{equation}
which are flags of length $d$.
\begin{thm}\label{Demazure-1}
For $w\in W^{\Lambda_0}$, we have
\begin{equation*}
\N_w=\bigsqcup_{\substack{\nu\in\cP \\
\ell(\nu)\leq d(w)}}SST_{\B_{<0}^\vee}(\nu^{\pi})_{\phi(w)} \times
SST_{\B_{>0}}(\nu^{\pi})_{\psi(w)}.
\end{equation*}
\end{thm}
\pf For convenience, let $\mathcal{S}_w$ be the right-hand side of
the above identity. We will use induction on $\ell(w)=|\lambda(w)|$.
When $\ell(w)=1$, i.e. $w=s_0$, it is clear. We assume that
$\ell(w)\geq 2$.

Choose $w'\in W^{\Lambda_0}$ such that $\ell(w')=\ell(w)-1$,
equivalently, $\lambda(w')\subset \lambda(w)$ with
$|\lambda(w)/\lambda(w')|=1$. Let $r$ be the residue of
$\lambda(w)/\lambda(w')$.
By (\ref{inductive relation}) and the induction hypothesis, we have
only to show that $\mathcal{S}_w=\bigcup_{k\geq
0}\tf_r^k\mathcal{S}_{w'} \setminus \{{\bf 0}\}$. We assume that
$\lambda(w)$ is as in (\ref{lambda w}) with $d=d(w)$.

\textsc{Case 1.} Suppose that $r=0$. Then we have $d\geq 2$ and
$\alpha_d=\beta_d=1$ and
\begin{equation*}
\lambda(w')=(\alpha_1,\ldots,\alpha_{d-1}|\beta_1,\ldots,\beta_{d-1}).
\end{equation*}

Let $(S,T)\in \mathcal{S}_{w'}$ be given, where ${\rm sh}(S)={\rm
sh}(T)=\eta^\pi$ for some $\eta\in\cP$ with $\ell(\eta)\leq
d(w')$. For $k\geq 1$, let $(S',T')=\tf_0^k(S,T)$ and
$\tau^\pi={\rm sh}(S')={\rm sh}(T')\in \cP^{\pi}$. By definition of
$\tf_0$, we have $\ell(\tau)\leq \ell(\eta)+1\leq d(w)$. If
$\ell(\tau)< d(w)$, then it is clear that $(S',T')\in
\mathcal{S}_{w'} \subset \mathcal{S}_{w}$. Assume that
$\ell(\tau)=d(w)$, that is, $\ell(\tau)=\ell(\eta)+1=d(w)$. Then the
first rows of $S'$ and $T'$ are filled only with $-1^\vee$ and $1$,
respectively, and the entries in the other rows of $S'$ and $T'$
still satisfy the flag conditions given by $\phi(w')$ and
$\psi(w')$, respectively. Since $\phi(w)=(1,\phi(w'))$ and
$\psi(w)=(1,\psi(w'))$, we have $(S',T')\in \mathcal{S}_{w}$.

Conversely, let $(S,T)\in \mathcal{S}_{w}$ be given with
$\tau^\pi={\rm sh}(S)={\rm sh}(T)\in \cP^{\pi}$. If
$\ell(\tau)<d(w)$, then $(S,T)\in \mathcal{S}_{w'}$. If
$\ell(\tau)=d(w)$, then the first rows of $S$ and $T$ are filled
only with $-1^\vee$ and $1$, respectively, and they are all
removable under successive application of $\te_0$. Hence, the
highest weight element $(S',T')$ in the $0$-string of $(S,T)$
belongs to $\mathcal{S}_{w'}$.\vskip 3mm

\textsc{Case 2.}  Suppose that $r\neq 0$. We may assume that $r>0$
since the argument for $r<0$ is almost the same. In this case, we
have $d(w')=d(w)$, and for some $1\leq i< d$
\begin{equation*}
\lambda(w')=(\alpha_1,\ldots,\alpha_i
,\ldots,\alpha_{d}|\beta_1,\ldots,\beta_i-1,\ldots,\beta_{d}),
\end{equation*}
where $\beta_i=r+1$.

Let $(S,T)\in \mathcal{S}_{w'}$ be given. Note that the entries of
$T$ in the $i$th row from the bottom are no more than $\beta_i-1=r$,
and no $r$ appears in the above rows. This implies that
$\tf_r^k(S,T)\in \mathcal{S}_{w}\cup\{{\bf 0}\}$ for $k\geq 0$.

Conversely, let $(S,T)\in \mathcal{S}_{w}$ be given. The entries of
$T$ in the $(i+1)$th row from the bottom are no more than
$\beta_{i+1}<\beta_i-1=r$. So any $r+1$ in the $i$th
row of $T$ from the bottom, if exists,  can be replaced by $r$ applying
$\te_r^k$ for some $k\geq 0$. This implies that $\te_r^k(S,T)\in
\mathcal{S}_{w'}$.
 \qed

\begin{cor} For $w\in W^{\Lambda_0}$, we have
\begin{equation*}
{\rm ch}\N_w=\sum_{\substack{\nu\in\cP \\
\ell(\nu)\leq
d(w)}}\widehat{s}_{\nu}(X_{\alpha(w)})\widehat{s}_{\nu}(Y_{\beta(w)}).
\end{equation*}
\end{cor}

Combining Theorem \ref{Demazure-0} and Theorem \ref{Demazure-1}, we
obtain a flagged version of the RSK correspondence and the Cauchy
identity.
\begin{cor}\label{main theorem} For $w\in W^{\Lambda_0}$, the map $\kappa$ in {\rm (\ref{kappa})}
gives a bijection
$$\{\,A\in\M\,|\ \lambda({\rm supp}(A))\subset
\lambda(w)\ \}\ \longrightarrow \bigsqcup_{\substack{\nu\in\cP \\
\ell(\nu)\leq d(w)}}SST_{\B_{<0}^\vee}(\nu^{\pi})_{\phi(w)} \times
SST_{\B_{>0}}(\nu^{\pi})_{\psi(w)},$$ when restricted to $\M_w$, and it commutes with $\te_i$ $(i\in\mathbb{Z})$. In
particular, we have
$$\sum_{\substack{S\subset\mathbb{N}^2 \\
w(S)\leq w}}\prod_{(i,j)\in S}\frac{x_iy_j}{(1-x_iy_j)}=\sum_{\substack{\nu\in\cP \\
\ell(\nu)\leq
d(w)}}\widehat{s}_{\nu}(X_{\alpha(w)})\widehat{s}_{\nu}(Y_{\beta(w)}).$$
\end{cor}

\begin{rem}{\rm (1) For $m,n\geq 1$, let $w_{m,n}$ be the element in $W^{\Lambda_0}$ such that
$\lambda(w_{m,n})=(n^m)$. In this case, we recover the usual RSK
correspondence with $m\times n$ matrices and the Cauchy identity
with variables $x_i,y_j$ ($1\leq i\leq m, 1\leq j\leq n$).


(2) For $S\subset \mathbb{N}^2$, let $r(S)$ be the diagonal length
of $\lambda(S)$, i.e. the length of $\alpha$ or $\beta$ when
$\lambda(S)=(\alpha|\beta)$. Then for $n\geq 1$, we have
\begin{equation*}
\sum_{\substack{S\subset\mathbb{N}^2 \\
r(S)\leq n}}\prod_{(i,j)\in S}\frac{x_iy_j}{(1-x_iy_j)}=\sum_{\substack{\nu\in\cP \\
\ell(\nu)\leq n}}{s}_{\nu}(X){s}_{\nu}(Y).
\end{equation*}
When multiplied by $e^{-n\Lambda_0}$ the right-hand side of the
identity is equal to the character of the irreducible highest weight
representation of $\widehat{\gl}_{\infty}$, a central extension of
$\gl_{\infty}$, with highest weight $-n\Lambda_0$, which is not
integrable \cite{KacR2}. Hence the left-hand side gives another
character formula for this highest weight module. Note that the
right-hand side has a Jacobi-Trudi type formula (see
\cite[Ex.7.16 d]{St} and \cite{K08} for its generalization to irreducible
$\widehat{\gl}_{\infty}$-modules with negative integral charges) and
a Weyl-Kac type formula \cite{K08-2}.

(3) From the correspondence between $\M_w$ and $\N_w$, we see that
the entries in the first columns of bitableux
 in $\N$ is determined only by the support of the coresponding matrix
in $\M$. This fact was also observed by Stanley
\cite[Ex.7.100]{St}  in a purely combinatorial way.  }
\end{rem}

\subsection{Demazure crystal $\B_w(n\Lambda_0)$}

\begin{prop}\label{Demazure-2} Let $w\in W^{\Lambda_0}$,
let $d=d(w)$, $\phi=\phi(w)$, and $\psi=\psi(w)$. Then for $n\geq 1$,
there is a bijection
\begin{equation*}
\B_w(n\Lambda_0) \longleftrightarrow \bigsqcup_{\substack{\nu\in\cP \\
\nu \subset (n^{d})}}SST_{\B_{<0}^\vee}(\nu^\pi)_{\phi} \times
SST_{\B_{>0}}(\nu^\pi)_{\psi}.
\end{equation*}
\end{prop}
\pf It follows from Remark \ref{imagePsi} and Theorem
\ref{Demazure-1}. \qed

\begin{rem}{\rm
Given $A=A(\bi,\bj)\in\M$ with $(\bi,\bj)\in\Omega$, let $c(A)$ be
the maximal length of decreasing subwords of $\bi$. It is well known
that $c(A)$ is equal to the number of columns in ${\bf P}(A)$ or
${\bf Q}(A)$ (cf.\cite{Kn}). By Remark \ref{imagePsi} and Theorem
\ref{Demazure-0}, the embedding $\Psi_n$ gives a bijection
\begin{equation*}
\B_w(n\Lambda_0) \longleftrightarrow \{\,A\,|\ \lambda({\rm
supp}(A))\subset \lambda(w),\ c(A)\leq n \ \}.
\end{equation*} }
\end{rem}

For $i\in\Z$ and $\lambda\in P$, let $D_{i}$ be the linear operator
on $\C[P]$ defined by
\begin{equation*}
D_i(e^{\lambda})=e^{\lambda}\cdot \frac{1-e^{-(1+\langle
\lambda,h_i\rangle)\alpha_i}}{1-e^{-\alpha_i}}.
\end{equation*}
The operators $D_i$ satisfy the braid relations, and hence for a reduced
expression of $w=s_{i_1}\cdots s_{i_r}\in W$, the operator
$D_w=D_{i_1}\cdots D_{i_r}$ is well-defined. By (\ref{Demazure}), we
have
\begin{equation*}
{\rm ch}\B_w(n\Lambda_0)=D_w(e^{n\Lambda_0}).
\end{equation*}
Combining with Proposition \ref{Demazure-2}, we obtain the following
combinatorial identity.

\begin{cor} Let $n, d\geq 1$, and
let $\alpha,\beta$ be two strict partitions of length $d$. Then
\begin{equation*}
{D_w(e^{n\Lambda_0})}e^{-n\Lambda_0}=\sum_{\substack{\nu\in\cP \\
\nu \subset
(n^{d})}}\widehat{s}_{\nu}(X_{\alpha})\widehat{s}_{\nu}(Y_{\beta}),
\end{equation*}
where $w$ is the unique element in $W^{\Lambda_0}$ such that
$\lambda(w)=(\alpha|\beta)$.
\end{cor}

\subsection{Crystals of symmetric matrices}\label{symmetric case}
From now on, let $\epsilon$ denote either $1$ or $2$. We put
\begin{equation}
\begin{split}
\widehat{\M}^\epsilon& =\{\,A\in\M\,|\,a_{-i^\vee,j}=a_{-j^\vee,i} \
\text{for $i,j\geq 1$},\ \text{$\epsilon$ divides $a_{-i^\vee,i}$
for $i\geq 1$} \,\}.
\end{split}
\end{equation}
For $i\in\mathbb{Z}_{\geq 0}$, let
\begin{equation}\label{folded operators}
\begin{split}
&\td{E}_0 =\te_0^\epsilon, \ \ \td{F}_0=\tf_0^\epsilon, \\
&\td{E}_i =\te_i\te_{-i},\ \ \td{F}_i =\tf_i\tf_{-i} \ \
(i\in\mathbb{Z}_{>0}).
\end{split}
\end{equation}
By similar arguments as in Proposition \ref{crystalM}, we can check
the following.
\begin{prop}\label{Mtilde}
\begin{itemize}
\item[(1)] $\widehat{\M}^{\epsilon}\cup\{{\bf 0}\}$ is invariant under
$\td{E}_i$ and $\td{F}_i$ for $i\in\mathbb{Z}_{\geq 0}$.

\item[(2)]
$\widehat{\M}^{\epsilon}=\{\,\td{F}_{i_1}\cdots\td{F}_{i_r}{\mathbb{O}}\,|\,r\geq
0,\  i_1,\ldots,i_r\in\mathbb{Z}_{\geq 0}\,\}\setminus\{{\bf 0}\}$.
\end{itemize}
\end{prop}

Put
\begin{equation*}
\begin{split}
\widehat{P}&=\{\,\lambda\in
P\,|\,\frac{1}{\epsilon}\langle\lambda,h_0\rangle\in\Z,\
\langle\lambda,h_i\rangle=\langle\lambda,h_{-i}\rangle\ (i\in\mathbb{Z}_{> 0})\,\},\\
\widehat{\Pi}&=\{\,\widehat{\alpha}_0=\epsilon\alpha_0, \ \
\widehat{\alpha}_i=\alpha_{i}+\alpha_{-i} \ \ (i\in\mathbb{Z}_{> 0})\,\}\subset \widehat{P}.\\
\end{split}
\end{equation*}
Then $\widehat{\Pi}$ is  a set of simple roots for the root system associated to
the affine Lie algebra $\mathfrak{b}_{\infty}$ (resp.
$\mathfrak{c}_{\infty}$) when $\epsilon=1$ (resp. $\epsilon=2$)
(cf.\cite{K}), where $\widehat{P}$ is a weight lattice. The
associated Dynkin diagrams are

\begin{center}
\hskip -1cm  \setlength{\unitlength}{0.16in}
\begin{picture}(24,4)
\put(2,2){\makebox(0,0)[c]{$\frak{b}_{\infty}\ :$}}

\put(5.6,2){\makebox(0,0)[c]{$\bigcirc$}}
\put(8,2){\makebox(0,0)[c]{$\bigcirc$}}
\put(10.4,2){\makebox(0,0)[c]{$\bigcirc$}}
\put(14.85,2){\makebox(0,0)[c]{$\bigcirc$}}
\put(17.25,2){\makebox(0,0)[c]{$\bigcirc$}}
\put(19.4,2){\makebox(0,0)[c]{$\bigcirc$}}
\put(8.35,2){\line(1,0){1.5}} \put(10.82,2){\line(1,0){0.8}}
\put(13.2,2){\line(1,0){1.2}} \put(15.28,2){\line(1,0){1.45}}
\put(17.7,2){\line(1,0){1.25}} \put(19.81,2){\line(1,0){0.9}}
\put(6.8,1.97){\makebox(0,0)[c]{$\Longleftarrow$}}
\put(12.5,1.95){\makebox(0,0)[c]{$\cdots$}}
\put(21.5,1.95){\makebox(0,0)[c]{$\cdots$}}
\put(5.6,1){\makebox(0,0)[c]{\tiny $\widehat{\alpha}_0$}}
\put(8,1){\makebox(0,0)[c]{\tiny $\widehat{\alpha}_1$}}
\put(10.4,1){\makebox(0,0)[c]{\tiny $\widehat{\alpha}_2$}}
\put(15,1){\makebox(0,0)[c]{\tiny $\widehat{\alpha}_{n-1}$}}
\put(17.4,1){\makebox(0,0)[c]{\tiny $\widehat{\alpha}_n$}}
\put(19.8,1){\makebox(0,0)[c]{\tiny $\widehat{\alpha}_{n+1}$}}
\end{picture}
\end{center}

\begin{center}
\hskip -1cm  \setlength{\unitlength}{0.16in}
\begin{picture}(24,4)
\put(2,2){\makebox(0,0)[c]{$\frak{c}_{\infty}\ :$}}

\put(5.6,2){\makebox(0,0)[c]{$\bigcirc$}}
\put(8,2){\makebox(0,0)[c]{$\bigcirc$}}
\put(10.4,2){\makebox(0,0)[c]{$\bigcirc$}}
\put(14.85,2){\makebox(0,0)[c]{$\bigcirc$}}
\put(17.25,2){\makebox(0,0)[c]{$\bigcirc$}}
\put(19.4,2){\makebox(0,0)[c]{$\bigcirc$}}
\put(8.35,2){\line(1,0){1.5}} \put(10.82,2){\line(1,0){0.8}}
\put(13.2,2){\line(1,0){1.2}} \put(15.28,2){\line(1,0){1.45}}
\put(17.7,2){\line(1,0){1.25}} \put(19.81,2){\line(1,0){0.9}}
\put(6.8,1.97){\makebox(0,0)[c]{$\Longrightarrow$}}
\put(12.5,1.95){\makebox(0,0)[c]{$\cdots$}}
\put(21.5,1.95){\makebox(0,0)[c]{$\cdots$}}
\put(5.6,1){\makebox(0,0)[c]{\tiny $\widehat{\alpha}_0$}}
\put(8,1){\makebox(0,0)[c]{\tiny $\widehat{\alpha}_1$}}
\put(10.4,1){\makebox(0,0)[c]{\tiny $\widehat{\alpha}_2$}}
\put(15,1){\makebox(0,0)[c]{\tiny $\widehat{\alpha}_{n-1}$}}
\put(17.4,1){\makebox(0,0)[c]{\tiny $\widehat{\alpha}_n$}}
\put(19.8,1){\makebox(0,0)[c]{\tiny $\widehat{\alpha}_{n+1}$}}
\end{picture}
\end{center}

For $i\in\mathbb{Z}_{\geq 0}$, let $\widehat{h}_i\in \widehat{P}^*$
be determined by $\langle \lambda, \widehat{h}_i\rangle = \langle
\lambda, h_i\rangle=\langle \lambda, h_{-i}\rangle$ if $i>0$, and
$\langle \lambda, \widehat{h}_0\rangle=\frac{1}{\epsilon}\langle
\lambda,h_0\rangle$ for   $\lambda\in \widehat{P}$. Then
$\widehat{\Pi}^\vee=\{\,\widehat{h}_i\,|\,i\in\mathbb{Z}_{\geq
0}\,\}$ is a set of simple coroots. As in Definition \ref{crystal
graph}, one may define an $x_{\infty}$-crystal
($x=\mathfrak{b,c}$) with respect to $\td{E}_i,\td{F}_i$,
$\widehat{\varepsilon}_i,\widehat{\varphi}_i$ ($i\in\mathbb{Z}_{\geq
0}$) and $\widehat{{\rm wt}}$. By Proposition \ref{Mtilde},
$\widehat{\M}^\epsilon$ is an $x_{\infty}$-crystal with highest
weight element $\mathbb{O}$. Here, for $A\in\widehat{\M}^\epsilon$
$\widehat{{\rm wt}}(A)={\rm wt}(A)\in \widehat{P}$,
$\widehat{\varepsilon}_i(A)=\varepsilon_i(A)$ and
$\widehat{\varphi}_i(A)=\varphi_i(A)$.

Since each $A$ in $\widehat{\M}^\epsilon$ is symmetric, we have
$\kappa(A)=(-S^\vee,S)$, where $S\in SST_{\B_{>0}}(\nu^\pi)$ for
some $\nu\in \cP$ with $\epsilon|\nu$, that is, $\epsilon|\nu_i$ for
$i\geq 1$ (cf.\cite{Kn,K07}), and $-S^\vee$ is the semistandard tableau
obtained by replacing each entry $i$ in $S$ with $-i^\vee$. Hence
the map $\widehat{\kappa} : A \mapsto S$ gives a bijection
\begin{equation}\label{rsk-1}
\widehat{\kappa}\ : \ \widehat{\M}^{\epsilon} \longrightarrow
\bigsqcup_{\substack{\nu\in \cP\\ \epsilon | \nu}}
SST_{\B_{>0}}(\nu^\pi ).
\end{equation}

\begin{prop}\label{Btilde} Put $\widehat{\Lambda}_0=\epsilon \Lambda_0$. For $n\geq 1$, let
\begin{equation*}
\widehat{\B}^{\epsilon}(n\widehat{\Lambda}_0)
=\{\,\td{F}_{i_1}\cdots\td{F}_{i_r}{H_{n\widehat{\Lambda}_0}}\,|\,r\geq
0, \ i_1,\ldots,i_r\in\mathbb{Z}_{\geq 0}\,\}\setminus\{{\bf 0}\}.
\end{equation*}
Then  $\Psi_{\epsilon
n}(\widehat{\B}^{\epsilon}(n\widehat{\Lambda}_0)\otimes
T_{-n\widehat{\Lambda}_0}) =\widehat{\M}^\epsilon \cap {\rm
Im}\Psi_{\epsilon n}$.
\end{prop}
\pf Let $\w\in \B(n\widehat{\Lambda}_0)$ be given with $\kappa\left(
\Psi_{\epsilon n}(\w\otimes
t_{-n\widehat{\Lambda}_0})\right)=(\w^\vee_{<0},\w_{>0})$. Then we
can check that
\begin{equation*}
\begin{split}
\w\in\widehat{\B}^\epsilon(n\widehat{\Lambda}_0) &
\Longleftrightarrow \w^\vee_{<0}=-\w^\vee_{>0} \Longleftrightarrow
\kappa^{-1}(\w^\vee_{<0},\w_{>0})\in\widehat{\M}^\epsilon,
\end{split}
\end{equation*}
where $-\w^\vee_{>0}$ is the tableau obtained from $\w_{>0}$ by
replacing each entry $i$ with $-i^\vee$. By Remark \ref{imagePsi},
we obtain the required identity. \qed

\begin{rem}\label{x-crystal}{\rm
  By Proposition \ref{Demazure-2}, (\ref{rsk-1}) and Proposition
\ref{Btilde}, we have a one-to-one correspondence
\begin{equation}\label{rsk-2}
\widehat{\B}^{\epsilon}(n\widehat{\Lambda}_0) \longleftrightarrow
\bigsqcup_{\substack{\nu\in \cP\\ \epsilon | \nu, \nu_1\leq
\epsilon n}} SST_{\B_{>0}}(\nu^\pi).
\end{equation}
By \cite[Theorem 5.1]{Kas96},
$\widehat{\B}^{\epsilon}(n\widehat{\Lambda}_0)$ is isomorphic to the
crystal graph of the integrable highest weight $x_{\infty}$-module
with highest weight $n\widehat{\Lambda}_0$, and Proposition
\ref{Btilde} implies that $\widehat{\M}^\epsilon$ is the limit of
$\widehat{\B}^{\epsilon}(n\widehat{\Lambda}_0)$.}
\end{rem}

Let $\sigma$ be the linear automorphism on $P$ defined by
$\sigma(\Lambda_0)=\Lambda_0$ and
$\sigma(\epsilon_i)=-\epsilon_{-i}$ for $i\in\Z$. Then
$\sigma(\alpha_i)=\alpha_{-i}$ for $i\in\Z$. Let
\begin{equation*}
\begin{split}
&\widehat{W}=\{\,w\in W\,|\,w\sigma=\sigma w\,\}. \\
\end{split}
\end{equation*}
Put $\widehat{s}_0=s_0$ and $\widehat{s}_i=s_is_{-i}$ for
$i\in\mathbb{Z}_{> 0}$. Then $\widehat{W}$ is the Coxeter group
generated by $\widehat{s}_i$ ($i\in\mathbb{Z}_{\geq 0}$), which is
isomorphic to the Weyl group of $x_{\infty}$ (see \cite[Section 5.2]{FSS}). Let $\widehat{W}^{\widehat{\Lambda}_0}$ be the set of minimal length
left coset representatives of $\widehat{W}_{\widehat{\Lambda}_0}$.
We have $\widehat{W}^{\widehat{\Lambda}_0}=\widehat{W}\cap
W^{\Lambda_0}$, and $\lambda(w)=\lambda(w)'$ or $\phi(w)=\psi(w)$
for $w\in \widehat{W}^{\widehat{\Lambda}_0}$.

For ${w}\in \widehat{W}^{\widehat{\Lambda}_0}$ with a reduced
expression ${w}=\widehat{s}_{i_1}\cdots\widehat{s}_{i_r}$
$(i_1,\ldots,i_r\in\mathbb{Z}_{\geq 0})$, let
\begin{equation}
\begin{split}
\widehat{\M}^{\epsilon}_{{w}}&=\{\,\td{F}_{i_1}^{m_1}\cdots\td{F}_{i_r}^{m_r}\mathbb{O}\,|\,m_1,\ldots,m_r\geq
0\,\}\setminus \{{\bf 0}\}, \\
\widehat{\B}^{\epsilon}_{{w}}(n\widehat{\Lambda}_0)&=\{\,\td{F}_{i_1}^{m_1}\cdots\td{F}_{i_r}^{m_r}H_{n\widehat{\Lambda}_0}\,|\,m_1,\ldots,m_r\geq
0\,\}\setminus \{{\bf 0}\}
\end{split}
\end{equation}
Since $\widehat{\B}^{\epsilon}(n\widehat{\Lambda}_0)$ is the crystal
graph of an integrable highest weight $x_{\infty}$-module (see
Remark \ref{x-crystal}), the Demazure crystal
$\widehat{\B}^{\epsilon}_{{w}}(n\widehat{\Lambda}_0)$ is
well-defined and so is $\widehat{\M}^{\epsilon}_{{w}}$ by
Proposition \ref{Btilde}.

\begin{prop}\label{Demazure-3} For ${w}\in
\widehat{W}^{\widehat{\Lambda}_0}$, we have
$$\widehat{\M}^{\epsilon}_{{w}}=\widehat{\M}^\epsilon\cap \M_w,\ \ \
\widehat{\kappa}(\widehat{\M}^{\epsilon}_{{w}})=\bigsqcup_{\substack{\nu\in \cP\\
\epsilon | \nu}} SST_{\B_{>0}}(\nu ^\pi)_{\phi({w})}.$$
\end{prop}
\pf Let us prove the first identity. Then the second one follows
from Proposition \ref{Demazure-1} and (\ref{rsk-1}). First, it is
clear by definition that
$\widehat{\M}^{\epsilon}_{{w}}\subset\widehat{\M}^\epsilon\cap
\M_w$. Suppose that $A\in \widehat{\M}^\epsilon\cap \M_w$ is given.
Let $w'=\widehat{s}_{i_2}\cdots\widehat{s}_{i_r}$, where
$w=\widehat{s}_{i_1}\cdots\widehat{s}_{i_r}$ is a reduced
expression. If $i_1>0$, then we have
$$\te_{i_1}A\neq
{\bf 0}\Longleftrightarrow \te_{-i_1}A\neq {\bf
0}\Longleftrightarrow \widetilde{E}_{i_1}A\neq {\bf 0},$$ since $A$
is symmetric. If $i_1=0$, then we have $\te_0A \neq {\bf 0}$ if and
only if $\widetilde{E}_0A\neq {\bf 0}$ by definition. From
(\ref{inductive relation}), it follows that
$\widetilde{E}^k_{i_1}A\in \M_{w'}$ for some $k\geq 0$. If we use
induction on $\ell(w)$, then we have $\widetilde{E}^k_{i_1}A\in
\widehat{\M}^\epsilon\cap \M_{w'}=\widehat{\M}^\epsilon_{w'}$, and
hence $A\in \widehat{\M}^\epsilon_{w}$. \qed

Put $Z=\{\,z_i=x_i y_i\,|\, i\geq 1\,\}$. For $S\subset \mathbb{N}^2$,
let $S'=\{\,(i,j)\,|\,(j,i)\in S\,\}$. Then Proposition \ref{Demazure-3} yields the
following identity.
\begin{cor} For ${w}\in
\widehat{W}^{\widehat{\Lambda}_0}$, we have
$$\sum_{\substack{S=S'\subset{\mathbb{N}^2} \\
w(S)\leq w}} \prod_{\substack{(i,i)\in
S}}\frac{z_i^\epsilon}{(1-z_i^\epsilon)} \prod_{\substack{(i,j)\in S
\\ i<j}}\frac{z_iz_j}{(1-z_iz_j)} =\sum_{\substack{\nu\in\cP \\
\epsilon|\nu \\ \ell(\nu)\leq
d(w)}}\widehat{s}_{\nu}(Z_{\alpha(w)}).$$
\end{cor}

\begin{cor}\label{Demazure-4}
For ${w}\in \widehat{W}^{\widehat{\Lambda}_0}$ and $n\geq 1$, there is in
one-to-one
correspondence
$$\widehat{\B}^{\epsilon}_{{w}}(n\widehat{\Lambda}_0)\ \ \ \longleftrightarrow \ \ \bigsqcup_{\substack{\nu\in \cP, \epsilon|\nu \\
\nu\subset (\epsilon n^{d({w})})}}
SST_{\B_{>0}}(\nu^\pi)_{\phi({w})}.$$
\end{cor}
\pf It follows from Proposition \ref{Demazure-2} and (\ref{rsk-2}).
\qed

For $i\in\mathbb{Z}_{\geq 0}$ and $\lambda\in \widehat{P}$, let
$\widehat{D}_{i}$ be the linear operator on $\C[\widehat{P}]\subset \C[{P}]$
defined by
\begin{equation*}
\widehat{D}_i(e^{\lambda})=e^{\lambda}\frac{1-e^{-(1+\langle
\lambda,\widehat{h}_i\rangle)\widehat{\alpha}_i}}{1-e^{-\widehat{\alpha}_i}}.
\end{equation*}
For a reduced expression of $w=\widehat{s}_{i_1}\cdots
\widehat{s}_{i_r}\in \widehat{W}$, we put
$\widehat{D}_w=\widehat{D}_{i_1}\cdots \widehat{D}_{i_r}$. Combining
the Demazure character formula ${\rm
ch}\widehat{\B}^{\epsilon}_{{w}}(n\widehat{\Lambda}_0)=\widehat{D}_{{w}}(e^{n\widehat{\Lambda}_0})$
with Corollary \ref{Demazure-4}, we obtain the following identity.
\begin{cor}
Let $n, d\geq 1$, and let $\alpha$ be a strict partition of length
$d$. Then
\begin{equation*}
{\widehat{D}_{{w}}(e^{n\widehat{\Lambda}_0})}\,{e^{-n\widehat{\Lambda}_0}}
=\sum_{\substack{\nu\in\cP, \epsilon|\nu \\
\nu \subset (\epsilon n^{d})}}\widehat{s}_{\nu}(Z_{\alpha}),
\end{equation*}
where ${w}$ is the unique element in
$\widehat{W}^{\widehat{\Lambda}_0}$ such that
$\lambda(w)=(\alpha|\alpha)$.
\end{cor}

\section{Plane partitions}

A {\it plane partition} is a collection of non-negative integers
$\pi=(\pi_{ij})_{i,j\geq 1}$ such that $\pi_{ij}\neq 0$ for only
finitely many $i,j$, and $\pi_{ij}\geq \pi_{i+1 j}$ and $\pi_{i
j}\geq \pi_{ij+1}$ for $i,j\geq 1$. A {\it shape of $\pi$} denoted
by ${\rm sh}(\pi)$ is a Young diagram determined by the support of
$\pi$, i.e. $\{\,(i,j)\,|\,\pi_{ij}\neq 0\,\}$. We may identify
$\pi$ with a tableau of shape ${\rm sh}(\pi)$ with entries in
$\mathbb{N}$ weakly decreasing in each row and column from left to
right and top to bottom, respectively.  Let $\mathcal{P}$ denote the set
of all plane partitions.

Let us recall the correspondence between $\M$ and $\mathcal{P}$
\cite{BK}. Let $A\in \M$ be given with $\kappa(A)=({\bf P}(A), {\bf
Q}(A))$. For each $k$, let $\lambda^{(k)}=(\alpha(k)|\beta(k))$ be a
partition where $\alpha(k)$ and $\beta(k)$ are strict partitions
given by reading the entries of the $k$th columns of ${\bf P}(A)$ and
${\bf Q}(A)$ from bottom to top (ignoring $-$ and $\vee$ in ${\bf
P}(A)$), respectively. Note that $\lambda^{(1)}\supset
\lambda^{(2)}\supset \cdots$. We define
$\pi(A)=(\pi(A)_{ij})_{i,j\geq 1}$ by $\pi(A)_{ij}=\left|
\{\,k\,|\,(i,j)\in\lambda^{(k)}\,\}  \right|$. It is easy to check
that $\pi(A)$ is a plane partition, and the mapping $A \mapsto
\pi(A)$ yields a bijection from $\M$ to $\mathcal{P}$. One may
identify a plane partition $\pi$ with a set of unit cubes, where
$\pi_{ij}$ cubes are stacked vertically at each position $(i,j)\in
{\rm sh}(\pi)$. Then $\lambda^{(k)}$ is the $k$th layer of $\pi(A)$
from the bottom.

For $\pi\in\mathcal{P}$, let $|\pi|=\sum_{i,j}\pi_{ij}$ and for each
$r\in\mathbb{Z}$, let ${\rm tr}_r(\pi)=\sum_{i\geq 1}\pi_{i i+r}$,
which is called the {\it $r$-trace of $\pi$} \cite{G,St73}. Let $q$
and $v_1,v_2,\ldots$ be formal variables. For a subset $X$ of
$\mathcal{P}$, the {\it norm} (resp. {\it trace}) {\it generating
function of $X$} is defined to be
\begin{equation*}\label{norm geneftn}
\begin{split}
\sum_{\pi\in X}q^{|\pi|}\ \ \ \text{and} \ \ \sum_{\pi\in
X}\prod_{r\in\mathbb{Z}}v_r^{{\rm tr}_r(\pi)},
\end{split}
\end{equation*}
respectively. Note that a norm generating function can be obtained
from a trace generating function by specializing $v_r=q$ for
$r\in\Z$.

For $n\geq 1$ and $\lambda\in\cP$, we put
\begin{equation}
\begin{split}
&\mathcal{P}(\lambda)=\{\,\pi\in\mathcal{P}\,|\,{\rm sh}(\pi)\subset
\lambda\,\},
\\
&\mathcal{P}_{\leq n}=\{\,\pi\in\mathcal{P}\,|\,\pi_{11}\leq n\,\},
\\
&\mathcal{P}(\lambda)_{\leq n}=\mathcal{P}_{\leq n}\cap
\mathcal{P}(\lambda).
\end{split}
\end{equation}
Note that $\mathcal{P}_{\leq 1}$ is the set of ordinary partitions
$\cP$, and hence a  $\gl_{\infty}$-crystal (cf.\cite{MM}), where for
$\lambda\in\cP$ and $r\in\mathbb{Z}$, $\tf_r\lambda$ is defined by a
partition $\mu$ such that $\mu/\lambda$ is a single box with residue
$r$, or ${\bf 0}$ if such $\mu$ does not exist.  If we assign the
weight of the empty partition to $0$, then it is easy to see that
$\mathcal{P}_{\leq 1}$ is isomorphic to $\B(\Lambda_0)\otimes
T_{-\Lambda_0}$. In general, we define a $\gl_{\infty}$-crystal
structure on $\mathcal{P}$ by identifying
$\pi=(\lambda^{(1)},\lambda^{(2)},\ldots,\lambda^{(n)})$
($\pi\in\mathcal{P}$) with an element of $\B(\Lambda_0)^{\otimes
n}\otimes T_{-n\Lambda_0}$ for a sufficiently large $n$, where
$\lambda^{(k)}$ is the $k$th layer of $\pi$.

\begin{prop}\label{PMiso} As $\gl_{\infty}$-crystals, we have
\begin{itemize}
\item[(1)] $\mathcal{P}\simeq \M$,

\item[(2)] $\mathcal{P}_{\leq n}
\simeq\B(n\Lambda_0)\otimes T_{-n\Lambda_0}$ for $n\geq 1$.
\end{itemize}
\end{prop}

\begin{cor}\label{trace generating function}\mbox{}
\begin{itemize}
\item[(1)] ${\rm ch}\M$ is the trace generating function of
$\mathcal{P}$.

\item[(2)] For $n\in\mathbb{N}$, $e^{-n\Lambda_0}{\rm
ch}\B(n\Lambda_0)$ is the trace generating function of
$\mathcal{P}_{\leq n}$.
\end{itemize}
\end{cor}
\pf For $\pi\in\mathcal{P}$, suppose that ${\rm
wt}(\pi)=-\sum_{r\in\mathbb{Z}}k_r\alpha_r$. Then $k_r$ is equal to
the $r$-trace of $\pi$. Identifying $e^{-\alpha_r}$ with $v_r$ in
${\rm ch}\M$ and ${\rm ch}\B(n\Lambda_0)$, we obtain the trace
generating functions for $\mathcal{P}$ and $\mathcal{P}_{\leq n}$,
respectively. \qed

\begin{rem}{\rm
By the celebrated Weyl-Kac character formula \cite{K}, the trace
generating function of $\mathcal{P}_{\leq n}$ is
$$\frac{\sum_{w\in
W}e^{w(n\Lambda_0+\rho)-n\Lambda_0-\rho}}{\prod_{\alpha\in\Delta^+}(1-e^{-\alpha})},$$
where $\rho\in\frak{h}^*$ is given by $\langle\rho,h_i\rangle=1$ for
all $i\in\mathbb{Z}$.
 The norm generating function for $\mathcal{P}$ and
$\mathcal{P}_{\leq n}$ are the corresponding principal
$q$-characters, which are obtained by putting $e^{-\alpha_i}=q$ for
$i\in\mathbb{Z}$. Then we recover
\begin{equation*}\label{nobound}
{\rm ch}_q\M=\frac{1}{\prod_{i\geq 1}(1-q^i)^i} \ \ \ \text{and} \ \
\ {\rm ch}_q\B(n\Lambda_0)=\frac{1}{\prod_{i\geq
1}(1-q^i)^{\min(i,n)}}
\end{equation*}
(see \cite{K} for evaluating $q$-characters), which are originally
due to MacMahon \cite{MacMahon}. }
\end{rem}

Now, the results in the previous section  give the following
representation theoretic interpretations on plane partitions with
bounded conditions.
\begin{prop}\label{planepartition} Let $\lambda\in\cP$ and
$n\geq 1$ be given. Let $w\in W^{\Lambda_0}$ be such that
$\lambda(w)=\lambda$. As $\gl_{\infty}$-crystals, we have
\begin{itemize}
\item[(1)] $\mathcal{P}(\lambda) \simeq \M_w$,

\item[(2)] $\mathcal{P}(\lambda)_{\leq n} \simeq \B_w(n\Lambda_0)\otimes T_{-n\Lambda_0}$.
\end{itemize}
\end{prop}
\pf Let $A\in\M$ be given. By Theorem \ref{Demazure-1}, we see that
$A\in\M_w$ if and only if $\lambda^{(1)}\subset \lambda(w)$. This
implies that $A\in\M_w$ if and only if ${\rm sh}(\pi(A))\subset
\lambda(w)$. The second isomorphism  follows from Proposition
\ref{Demazure-2}. \qed\vskip 3mm

By Corollary \ref{trace generating function} and Proposition
\ref{planepartition}, we obtain the generating functions for plane
partitions bounded by a given shape as follows.
\begin{cor}\label{generating functions}
Let $\lambda\in\cP$ and $n\geq 1$ be given. Let $w\in W^{\Lambda_0}$
be such that $\lambda(w)=\lambda$.
\begin{itemize}
\item[(1)] The trace generating function of
$\mathcal{P}(\lambda)$ is
\begin{equation*}
\sum_{\substack{S\subset\Delta(\mathfrak{u}_-) \\
\lambda(S)\subset \lambda}}\prod_{\alpha\in
S}\frac{e^\alpha}{(1-e^\alpha)}.
\end{equation*}

\item[(2)] The trace generating function of
$\mathcal{P}(\lambda)_{\leq n}$ is
$e^{-n\Lambda_0}D_w(e^{n\Lambda_0})$.
\end{itemize}
\end{cor}

\begin{rem}{\rm

(1) There are determinantal formulas for the norm and trace
generating functions of various classes of plane partitions
 including
$\mathcal{P}(\lambda)$ and $\mathcal{P}(\lambda)_{\leq n}$ (see
\cite{Kr} for a most general form and the references therein for the
previous works by other people). Also there are evaluations of those
determinants into nice product forms for some special classes of
plane partitions. But there is no such formula for
$\mathcal{P}(\lambda)$ and $\mathcal{P}(\lambda)_{\leq n}$ as far as
we know.

(2) A representation theoretic approach to plane partitions was
first introduced by Proctor \cite{Pr}, where the norm generating
functions for $\mathcal{P}((u^v))_{\leq n}$ was proved to be the
$q$-dimension of the irreducible $\frak{sl}_{u+v}$-module with
highest weight $n\omega$ ($\omega$ is the $u$th fundamental weight).
}
\end{rem}

A plane partition $\pi=(\pi_{ij})$ is called {\it symmetric} if
$\pi_{ij}=\pi_{ji}$ for all $i,j\geq 1$. Similarly for
$\epsilon=1,2$, $n\geq 1$ and $\lambda\in\cP$ with
$\lambda=\lambda'$, we put
\begin{equation}
\begin{split}
& \widehat{\mathcal{P}}^\epsilon=\{\,\pi\in\mathcal{P}\,|\,\text{$\pi$ is
symmetric   and $\epsilon$ divides $\pi_{ii}$ for
all $i\geq 1$} \,\}, \\
&\widehat{\mathcal{P}}(\lambda)^\epsilon=\mathcal{P}(\lambda)\cap
\widehat{\mathcal{P}}^\epsilon,
\\
&\widehat{\mathcal{P}}^{\epsilon}_{\leq n}=\mathcal{P}_{\leq n}\cap
\widehat{\mathcal{P}}^\epsilon,
\\
&\widehat{\mathcal{P}}(\lambda)^{\epsilon}_{\leq
n}=\mathcal{P}_{\leq n}(\lambda)\cap\widehat{\mathcal{P}}^\epsilon.
\end{split}
\end{equation}
It is clear that $\widehat{\M}^\epsilon$ is in one-to-one
correspondence with $\widehat{\mathcal{P}}^\epsilon$. We define for
a subset $X$ of $\widehat{\mathcal{P}}^\epsilon$ the norm (resp.
trace) generating function of $X$ by
\begin{equation*}
\begin{split}
\sum_{\pi\in X}q^{|\pi|}\ \ \ \text{and} \ \ \sum_{\pi\in
X}\prod_{r\geq 0}v_r^{{\rm tr}'_r(\pi)},
\end{split}
\end{equation*}
respectively, where ${\rm tr}'_r(\pi)={\rm tr}_r(\pi)$ for $r\geq 1$
and ${\rm tr}'_0(\pi)=\epsilon^{-1}{\rm tr}_0(\pi)$.

As in Section \ref{symmetric case}, we assume that $x=\frak{b}$ if
$\epsilon=1$ and $x=\frak{c}$ if $\epsilon=2$. We define an
$x_{\infty}$-crystal structure on
$\widehat{\mathcal{P}}^\epsilon\subset \mathcal{P}$ with
$\widetilde{E}_i$, $\widetilde{F}_i$ for $i\geq 0$ (\ref{folded
operators}). Similar to Proposition \ref{PMiso} and Proposition
\ref{planepartition}, we can prove the following.

\begin{prop}\label{symmetric plane partitions} Let $n\geq 1$ and $\lambda\in\cP$ with $\lambda=\lambda'$ be given,
and let $w\in\widehat{W}^{\widehat{\Lambda}_0}$ be such that
$\lambda(w)=\lambda$. As $x_{\infty}$-crystals, we have
\begin{itemize}
\item[(1)] $\widehat{\mathcal{P}}^{\epsilon}\simeq \widehat{\M}^\epsilon$,

\item[(2)] $\widehat{\mathcal{P}}^{\epsilon}(\lambda)\simeq
\widehat{\M}^\epsilon_w$,

\item[(3)] $\widehat{\mathcal{P}}^{\epsilon}_{\leq n}
\simeq\widehat{\B}^\epsilon(n\widehat{\Lambda}_0)\otimes
T_{-n\widehat{\Lambda}_0}$,

\item[(4)] $\widehat{\mathcal{P}}(\lambda)^{\epsilon}_{\leq n}
\simeq\widehat{\B}_w^\epsilon(n\widehat{\Lambda}_0)\otimes
T_{-n\widehat{\Lambda}_0}$.
\end{itemize}
\end{prop}

\begin{rem}{\rm
(1) For $\pi\in\widehat{\mathcal{P}}^\epsilon$, we have ${\rm
wt}(\pi)=-\sum_{r\geq 0}k_r\widehat{\alpha}_r$, where $k_r={\rm
tr}'_r(\pi)$. Hence, identifying $e^{-\widehat{\alpha}_r}$ with
$v_r$, we obtain the trace generating functions for symmetric plane
partitions in Proposition \ref{symmetric plane partitions} as the
characters of the corresponding $x_{\infty}$-crystals. The norm
generating functions can be obtained by specializing
$e^{-\widehat{\alpha}_r}=q^2$ for $r\geq 1$ and
$e^{-\widehat{\alpha}_0}=q^\epsilon$.

(2) The norm generating function of
$\widehat{\mathcal{P}}(m^m)^{1}_{\leq n}$ was conjectured by
MacMahon \cite{MacMahon}, and it was proved by Andrews \cite{A} and
Macdonald \cite{Mac95}. It was observed by Proctor \cite{Pr,Pr90}
that the norm generating function of
$\widehat{\mathcal{P}}(m^m)^{\epsilon}_{\leq n}$ is the
$q$-dimension of irreducible representation of the complex simple
Lie algebra $\frak{so}(2m+1)\subset \frak{b}_{\infty}$ or
$\frak{sp}(2m)\subset \frak{c}_{\infty}$ with highest weight
corresponding to $n\widehat{\Lambda}_0$ (following our notation).

(3) Recently Tingley \cite{T} gave a nice representation theoretic
interpretation of cylindric plane partitions in terms of crystal
graphs for affine Lie algebra $\widehat{\frak{sl}}_n$ and its generating function. It would be
interesting to find an application of affine Demazure crystals to
cylindric plane partitions.}
\end{rem}


{\small
\begin{thebibliography}{BSS}
\bibitem{A}
G. E. Andrews, {\em Plane partitions. I. The MacMahon conjecture},
Studies in foundations and combinatorics, pp. 131--150, Adv. in
Math. Suppl. Stud., 1, Academic Press, New York-London, 1978.

\bibitem{BK}
E.~A. Bender, D.~E. Knuth, {\em Enumeration of plane partitions}, J.
Combinatorial Theory Ser. A \textbf{13} (1972), 40--54.


\bibitem{Br}
J. Brundan, {\em Dual canonical bases and Kazhdan-Lusztig
polynomials}, J. Algebra \textbf{306} (2006), 17--46.

\bibitem{DK}
V.~I. Danilov, G.~A. Koshevoy, {\em Bi-crystals and crystal
$(GL(V),GL(W))$ duality}, RIMS preprint, (2004) no. 1458.


\bibitem{FSS}
J. Fuchs, A. N. Schellekens, C. Schweigert, {\em From Dynkin diagram
symmetries to fixed point structures}, Comm. Math. Phys.
\textbf{180} (1996), 39--97.

\bibitem{G}
E. R. Gansner, {\em  The Hillman-Grassl correspondence and the
enumeration of reverse plane partitions}, J. Combin. Theory Ser. A
\textbf{30} (1981), 71--89.

\bibitem{H} R.~Howe, {\em Remarks on classical invariant theory}, Trans.
AMS {\bf 313} (1989), 539--570.

\bibitem{K} V. G.~Kac, {\em Infinite-dimensional Lie algebras}, Third
edition, Cambridge University Press, Cambridge, 1990.

\bibitem{KacR2}
V. G. Kac, A. Radul, {\em Representation theory of the vertex
algebra $W\sb {1+\infty}$}, Transform. Groups \textbf{1} (1996),
41--70.

\bibitem{Kas90}
M. Kashiwara, {\em Crystalizing the $q$-analogue of universal
enveloping algebras}, Comm. Math. Phys. \textbf{133} (1990),
249--260.

\bibitem{Kas93}
M. Kashiwara, \emph{Crystal bases and Littelmann's refined Demazure
character formula}, Duke Math. J. \textbf{71} (1993),  839--858.

\bibitem{Kas94}
M. Kashiwara, \emph{On crystal bases}, Representations of groups,
CMS Conf. Proc., 16, Amer. Math. Soc., Providence, RI, (1995),
155--197.

\bibitem{Kas96}
M. Kashiwara, {\em Similarity of crystal bases}, Contemp. Math.
\textbf{194} (1996), 177--186.

\bibitem{KN}
M. Kashiwara, T. Nakashima, {\em Crystal graphs for representations
of the $q$-analogue of classical Lie algebras}, J. Algebra
\textbf{165} (1994), 295--345.

\bibitem{Kn}
D. E. Knuth, {\em Permutations, matrices, and  generalized Young
tableaux}, Pacific J. Math. \textbf{34} (1970),  709--727.

\bibitem{Kr}
C. Krattenthaler, {\em Generating functions for plane partitions of
a given shape}, Manuscripta Math. \textbf{69} (1990), 173--201.

\bibitem{K07}
J.-H. Kwon, {\em Crystal graphs for Lie superalgebras and Cauchy
decomposition}, J. Algebraic Combin. \textbf{25} (2007), 57--100.

\bibitem{K08}
J.-H. Kwon, {\em Rational semistandard tableaux and character
formula for the Lie superalgebra
$\widehat{\frak{gl}}_{\infty|\infty}$}, Adv. Math. \textbf{217}
(2008), 713-739.

\bibitem{K08-2}
J.-H. Kwon,  {\em A combinatorial proof of a Weyl-type formula for
hook Schur polynomials},  J. Algebraic Combin. \textbf{28} (2008), 439--459

\bibitem{La}
A. Lascoux, {\em Double crystal graphs}, Studies in Memory of Issai
Schur, Progress in Math. \textbf{210}, Birkha\"{u}ser (2003),
95--114.

\bibitem{Lc}
C. Lecouvey, {\em  Crystal bases and combinatorics of infinite rank
quantum groups}, preprint, arXiv:math/0604636.

\bibitem{Ln}
C. Lenart, {\em A unified approach to combinatorial formulas for
Schubert polynomials}, J. Algebraic Combin. \textbf{20} (2004),
263--299.

\bibitem{vL}
M. A. A. van Leeuwen,  {\em Double crystals of binary and integral
matrices}, Electron. J. Combin. \textbf{13} (2006).

\bibitem{Mac}
I. G. Macdonald, {\em Notes on Schubert polynomilas}, Laboratoire de
combinatoire et d'informatique math\'{e}matique (LACIM),
Universit\'{e} du Qu\'{e}bec \`{a} Montr\`{e}al, Montreal, 1991.

\bibitem{Mac95}
I. G. Macdonald, {\em Symmetric functions and Hall polynomials},
Oxford University Press, 2nd ed.,  1995.

\bibitem{MacMahon}
P. A. MacMahon, {\em Combinatory analysis}, vols. 1 and 2, Cambridge
University Press, Canbridge, 1915, 1916;reprinted in one volume by
Chelsea, New York, 1960.

\bibitem{MM}
K.~Misra, T.~Miwa, \emph{Crystal base for the basic representation
of ${U}_q(\widehat{\mathfrak{sl}}(n))$}, Comm. Math. Phys.
\textbf{134} (1990),  79--88.

\bibitem{Pr}
R. A. Proctor, {\em Bruhat lattices, plane partition generating
functions, and minuscule representations}, European J. Combin. \textbf{5}
(1984),   331--350.

\bibitem{Pr90}
R. A. Proctor, {\em New symmetric plane partition identities from
invariant theory work of De Concini and Procesi}, European J.
Combin. \textbf{11} (1990), 289--300.


\bibitem{St73}
R. P. Stanley, {\em The conjugate trace and trace of a plane
partition}, J. Combinatorial Theory Ser. A \textbf{14} (1973),
53--65.

\bibitem{St}
R. P. Stanley, {\em Enumerative Combinatorics}, vol. 2, Cambridge
University Press, 1998.



\bibitem{SS}
B. E. Sagan, R. Stanley, {\em Robinson-Schensted algorithms for skew
tableaux}, J. Combin. Theory Ser. A \textbf{55} (1990), 161--193.

\bibitem{T}
P. Tingley, {\em   Three combinatorial models for affine
${\frak{sl}}_n$ crystals, with applications to cylindric plane
partitions}, Int.~Math.~Res.~Not. \textbf{143} (2007).
\end{thebibliography}}
\end{document}